\UseRawInputEncoding
\documentclass[a4paper,10pt, reqno]{amsart}
\usepackage{amssymb,latexsym,amsmath,epsfig,amsthm} %% Add other packages as necessary 

\theoremstyle{plain}
\newtheorem{theorem}{Theorem}
\newtheorem*{theorem*}{Theorem}
\newtheorem{corollary}{Corollary}
\newtheorem*{corollary*}{Corollary}
\newtheorem{lemma}{Lemma}
\newtheorem*{lemma*}{Lemma}

\newtheorem*{proposition*}{Proposition}

\newtheorem*{conjecture*}{Conjecture}
\theoremstyle{definition}
\newtheorem{definition}{Definition}
\newtheorem*{definition*}{Definition}
\theoremstyle{remark}
\newtheorem{remark}{Remark}
\newtheorem*{remark*}{Remark}

\newtheorem*{property*}{Property}
\newtheorem{example}{Example}
\newtheorem*{example*}{Example}

\begin{document}
\title[Fractal sets]{More on one class of fractals ({Some fractal properties of sets
having the Moran structure})}%
\author{Symon Serbenyuk}
\address{
  45~Shchukina St. \\
  Vinnytsia \\
  21012 \\
  Ukraine}
\email{simon6@ukr.net}

\subjclass[2010]{28A80, 11K55, 26A09}

% Key words
\keywords{ Fractal,  Cantor-like set, Moran structure, Hausdorff  dimension, self-similar set, s-adic representation, nega-s-adic representation, alternating Cantor series, mixed s-adic series,  nega-s-adic Cantor series.}

\begin{abstract}
This article is devoted to  sets having the Moran structure. The main attention is given to topological, metric, and fractal properties of certain sets whose  elements have restrictions on using  digits or combinations of digits in own representations.

\end{abstract}
\maketitle

\section{Introduction}

In 1977, the notion ``fractal" was considered by B.~Mandelbrot in  \cite{Mandelbrot1977}. A fractal in the wide sense is a set whose topological dimension does not coincide with the Hausdorff  dimension (the fractal dimension), and in the narrow sense it is a set that  has the fractional fractal dimension. 

Fractals are the most appropriate mathematical models of natural objects. 
The importance of fractals lies in modeling of physical and biological processes, and  also fractal is a strictly mathematical notion that unites various mathematical objects, e.g. continuous nowhere differentiable  functions, singular distributions, curves and surfaces that do not have the tangent at any point, etc. (see \cite{{Falconer1997}, {Falconer2004}, {Mandelbrot1977}}). 
Indeed,  the following examples are natural examples of fractals: the rings around planets (such fractals have the property of self-similarity),   the snow cover in a mountain region, linear  lightning, cloud borders, forms of coast lines or rivers. In fact, one can model coast lines and rivers by continuous nowhere differentiable functions. One of the oldest mathematical examples of fractals is the Cantor set
$$
C=\left\{x: x=\sum^{\infty} _{n=1}{\frac{\alpha_n}{3^n}}, \alpha_n\in\{0,2\}\right\}.
$$ 
This set was intoduced by G.~Cantor in 1883. The last fractal is the part of such well-known fractals as the Sierpinski  carpet, the Koch snowflake.

Fractal sets are widely applicated  in computer design, algorithms of the compression to information, quantum mechanics, solid-state physics, analysis and categorizations of signals of  various forms appearing in different areas (e.g. the analysis of exchange rate fluctuations in economics),  etc. However, for much classes of fractals the problem of the Hausdorff dimension calculation is difficult and the estimate of parameters  on which the Hausdorff  dimension of certain classes of fractal sets depends is left out of consideration. 

The aim of this survey is to give some local and fractal properties of certain Moran sets. One can note that Moran sets play an important role in multifractal analysis/formalism and especially the refined multifractal formalism (for example, see the papers \cite{2021-1, Wu2005, Wu2005(2), XW2008} and references therein).

The multifractal analysis is a natural framework to finely describe geometrically the heterogeneity in the distribution at small scales of the measures on a metric space. The multifractal analysis was proved to be a very useful technique in the analysis of
measures, both in theory and applications.  Also, the multifractal and the fractal analysis allows one to perform a certain classification of singular measures. One can note that it was proved that singular distributions of probabilities are dominant for many classes of random variables. Possible applications in the spectral theory of self-adjoint operators serve as an additional stimulus for a further investigation of singularly continuous measures \cite{DS2020-ERA}. For example, one can note the following researches of singular measures: singularity of Hewitt--Stromberg measures on Bedford--McMullen carpets \cite{AB2021-AM}, the mutual singularity of certain measures (see \cite{{Yuan2019-N}, {DS2020-ERA}, {XW2008-Fr}, {Selmi2021-BPASciM}, {DS2021-RM}} and references therein ), dimensions of measures \cite{{RS2021-JGA}, {HLW2020-JMAA}, {Selmi2022-ASM}}.

Olsen \cite{Olsen1995} introduced a general form of multifractal formalism, to interpret the statistical scaling properties of singular measures where the total mass or energy is spread over regions of phase space in an irregular way. The multifractal formalism aims at expressing the dimensions (the Hausdorff and packing dimensions) of the level sets in terms of the Legendre transform of some free energy function in analogy with the usual thermodynamic theory (\cite{{Selmi2022-ASM}, {Yuan2019-N}, {2021-1}} and references therein ).

The multifractal formalism has been proved rigorously for random and non-random self-similar measures, for self-conformal measures, for self-affine and forMoran measures (see \cite{{2021-1}, {DS2021-Pr}} and references therein). Certan researches  are devoted to new multifractal formalism for which the classical multifractal formalism does not hold. For example, the paper \cite{2021-1}] deals with a multifractal formalism based on the Hewitt--Stromberg measures and that this formalism is completely parallel to Olsen's multifractal formalism which is based on the Hausdorff and packing measures. 

Among fractal geometrical objects, Moran's types play an important role in explaining many situations, in pure mathematics as the general context of Cantor's, and in applied physics as a suitable context for studying scaling laws. These sets may be understood as attractors for dynamical systems, electrical circuits, and also smart cities where fractals are nowadays sophisticated tools in their modeling.
Fractals such as Cantor, and generally Moran's types are also applied in understanding physical properties at different molecular levels, such as nonmaterial composites, crystal growth, and structure, porous materials, etc. \cite{DSM2021-EPhJ}. Finally, one can note some investigations in multifractal analysis of Moran sets: multifractal properties  of homogeneous Moran fractals associated with Fibonacci sequence \cite{Wu2005(2)}, multifractal properties \cite{{Wu2005}, {WX2011-ChSF}}.

Consider space $\mathbb R^n$. In \cite{Moran1946}, P. A. P. Moran introduced the following construction of sets and calculated the Hausdorff dimension of the limit set 
\begin{equation}
\label{eq: Cantor-like set}
E=\bigcap^{\infty} _{n=1}{\bigcup_{i_1,\dots , i_n\in A_{0,p}}{\Delta_{i_1i_2\ldots i_n}}}. 
\end{equation}
 Here $p$ is a fixed positive integer, $A_{0,p}=\{1, 2, \dots , p\}$, and sets $\Delta_{i_1i_2\ldots i_n}$ are basic sets  having  the following properties:
\begin{itemize}
\item any set $\Delta_{i_1i_2\ldots i_n}$ is closed and disjoint;
\item for any $i\in A_{0,p}$ the condition $\Delta_{i_1i_2\ldots i_ni}\subset\Delta_{i_1i_2\ldots i_n}$ holds;
\item 
$$
\lim_{n\to\infty}{d\left(\Delta_{i_1i_2\ldots i_n}\right)}=0, \text{where $d(\cdot)$ is the diameter of a set};
$$
\item each basic set is the closure of its interior;
\item at each level the basic sets do not overlap (their interiors are disjoint);
\item any basic set $\Delta_{i_1i_2\ldots i_ni}$ is geometrically similar to $\Delta_{i_1i_2\ldots i_n}$;
\item
$$
\frac{d\left(\Delta_{i_1i_2\ldots i_ni}\right)}{d\left(\Delta_{i_1i_2\ldots i_n}\right)}=\sigma_i,
$$
where $\sigma_i\in (0,1)$ for $i=\overline{1,p}$.
\end{itemize}

The Hausdorff-Besicovitch dimension $\alpha_0$ of the set $E$ is the unique root of the following equation
$$
\sum^{p} _{i=1}{\sigma^{\alpha_0} _i}=1.
$$

It is easy to see that set \eqref{eq: Cantor-like set} is a Cantor-like set and a self-similar fractal. The set $E$  is called \emph{the Moran set}.

Let us consider   the second definition of the Moran set given by  Hua et al.~(\cite{HRW2000}).
\begin{definition*}{(Definition of  Hua et al.).}[ {\cite{Serbenyuk2021}}]
Let $(n_k)$ be a sequence of positive integers, $J\in\mathbb R^n$ be a compact set with nonempty interior, and $(\Phi_k)$ be a sequence of positive real vectors with
$\Phi_k=(\sigma_{k,1},\sigma_{k,2},\dots , \sigma_{k,n_k}),$
where $k\in\mathbb N$ and
$$
\sum^{n_k} _{j=1}{\sigma_{k,j}}<1.
$$
A  set  of the form
$$
E=\bigcap^{\infty} _{k=1}{\bigcup_{i_1,\dots , i_k\in A_{0,n_k}}{\Delta_{i_1i_2\ldots i_k}}},
$$
where $A_{0,n_k}=\{1, 2, \dots , n_k\}$, is called \emph{the Moran
set associated with the collection $F$}. 
Here
$$
F=\bigcup^{\infty} _{k=0}{F_k}=\bigcup^{\infty} _{k=0}{\left\{J_{\sigma}:=\Delta_{i_1i_2...i_k}: k\in\mathbb N, i_k\in\{1, 2, \dots , n_k\}      \right\}}
$$
\emph{The collection $F$ fulfills the Moran structure} provided it satisfies the following Moran Structure Conditions (MSC):

\begin{enumerate}
\item 
\label{property-1}
$J_{\varnothing}=J$.
\item
\label{property-2}
An arbitrary $J_{\sigma}$ is geometrically similar to $J$.
\item
\label{property-3}
 For any $i,j\in \{1,2, \dots , n_{k+1}\}$ such  that $i\ne j$, the conditions 
$$
\Delta_{i_1i_2\ldots i_ki}\subset\Delta_{i_1i_2\ldots i_k}, ~~~\Delta_{i_1i_2\ldots i_ki}\cap\Delta_{i_1i_2\ldots i_kj}=\varnothing
$$
hold.
\item
\label{property-4} For any $j\in \{1,2, \dots , n_{k+1}\}$,
$$
\frac{d\left(\Delta_{i_1i_2\ldots i_kj}\right)}{d\left(\Delta_{i_1i_2\ldots i_k}\right)}=\sigma_{k+1, j}.
$$
\end{enumerate}
The elements of $F_k$ are called \emph{the basic elements
of order k of the Moran set $E$}, and the elements of $F$ are called \emph{ the basic elements of the
Moran set $E$}.
\end{definition*}

\begin{remark}
Let us note that the main difference between definitions of Moran  and Hua is Property~\ref{property-4} in MSC.
\end{remark}

Let $M=M(J,(n_k), (\Phi_k))$ be a class of Moran sets satisfying MSC~\ref{property-1}--\ref{property-4}. It is known that one can define a sequence $(\alpha_k)$, where $\alpha_k$ satisfies the equation
$$
\prod^{k} _{i=1}{\sum^{n_i} _{j=1}{\sigma^{\alpha_k} _{i,j}}}=1.
$$
Also, suppose that 
$$
\alpha_{*}=\liminf_{k\to\infty}{\alpha_k}, ~~~ \alpha^{*}=\limsup_{k\to\infty}{\alpha_k};
$$
$$
c_{*}=\inf_{i,j}{\sigma_{i,j}}, ~~~~c^{*}=\sup_{i,j}{\sigma_{i,j}}.
$$

Much research has been devoted to Moran-like constructions and Cantor-like sets (for example, see \cite{PW1996, HRW2000, PS1995, Mance2015, DU2014, DU2014(2)} and references therein). For example, in \cite{PS1995},  the one parameter family of Cantor sets
$$
\Lambda(\lambda)=\left\{x: x=\sum^{\infty} _{k=1}{i_k\lambda^k}, i_k\in S\subset \{0,1,\dots , s-1\}, s\in \mathbb N~ \text{is a fixed number}\right\}
$$
 is investigated.
\begin{theorem}[\cite{PS1995}]
Suppose that the condition $s-1<(l-1)^2$ holds. Here $l$ is the cardinality of the set $S=\{s_1,\dots , s_l\}$, i.e., $l=|S|$. Then for almost all $\lambda\in\left[\frac{1}{s},\frac{1}{l}\right]$ (with respect to Lebesgue measure) we have that 
$$
\alpha_0 \left(\Lambda(\lambda)\right)=\frac{\log l}{-\log \lambda}.
$$
\end{theorem}
It is easy to see that we obtain the case of  classical s-adic representation \eqref{eq: s-adic expansion} whenever $\lambda=\frac{1}{s}$. In this case, we get
$$
\alpha_0 \left(\Lambda(\lambda)\right)=\log_{s}{l}.
$$
The following theorem generalizes the last result.

Let $D=(d_n)$ be a fixed sequence of positive integers such that $d_n>1$ for all $n\in\mathbb N$, $\varepsilon_n\in A_{d_n}=\{0,1,\dots , d_n-1\}$. Series of the form
$$
\sum^{\infty} _{n=1}{\frac{\varepsilon_n}{d_1d_2\cdots d_n}}
$$
are  \emph{ Cantor series} introduced  by G. Cantor in \cite{Cantor1}. These series are  generalizations of  s-adic expansion \eqref{eq: s-adic expansion}, i.e., a Cantor series is the s-adic expansion whenever $d_n=const=s$ for all $n\in\mathbb N$.  

\begin{theorem}[\cite{Mance2015}]
Suppose that $D=(d_n)$ be a fixed sequence of positive integer numbers $d_n>1$, $\lim_{n\to\infty}{\frac{\log{d_n}}{\log{d_1\cdots d_n}}}=0$, $I_j\subseteq\{0,1, \dots , d_j-1\}$, $\mathcal{I}=(I_n)$. Then
$$
\alpha_0\left(\mathcal{R}_{\mathcal{I}}(D)\right)=\alpha_0\left(\left\{x: x=\sum^{\infty} _{n=1}{\frac{\varepsilon_n}{d_1d_2\cdots d_n}}, \varepsilon_n\in I_n\right\}\right)={\lim\inf }_{n\to\infty}\frac{\log{\prod^{n} _{j=1}{|I_j|}}}{\log{\prod^{n} _{j=1}{d_j}}}.
$$
\end{theorem}

The present survey  is devoted to fractal sets, whose elemens defined by expansions related with some cases of positive and alternating Cantor series and their images under the action of  certain singular distributions. The main attention is given to topological and metric properties of these sets, and also parameters under which depends  the Hausdorff  dimension of such sets.  Sets considered in this paper, are determined by certain restrictions on using combinations of digits in representations of them elements.
Also, the main attention is given to results obtained in the parers \cite{Symon1, Symon2, S. Serbenyuk 2013} published into Ukrainian. 

Let us remark that, in September 2011 and February 2012,  results of \cite{Symon1, Symon2} were presented by the author in the reports 
``The main topological, metric properties of one set of numbers such that it is defined by the s-adic representation with restrictions" and  ``The main topological, metric properties of one set defined by the nega-s-adic  and s-adic representation with a parameter,  and using this set"  at the fractal analysis  seminar of the Institute of Mathematics of NAS of Ukraine and the National Pedagogical Dragomanov University (archive of reports is available here:\\ http://www.imath.kiev.ua/events/index.php?seminarId=21\&archiv=1). In 2012, results of the papers \cite{Symon1, Symon2} were presented in the conference abstracts \cite{S. Serbenyuk abstract 2, S. Serbenyuk abstract 3, S. Serbenyuk abstract 5}. Also, the main results of these papers were published into English as the preprint~\cite{{S. Serbenyuk 2017  fractals}}. 

\section{Definitions}
We begin with definitions of several representations of real numbers and certain series. 

Let $1<s$ be a fixed positive integer, $A=\{0, 1,\dots , s-1\}$ be an alphabet of the s-adic or nega-s-adic numeral system, and $A_0=A \setminus \{0\}=\{1,2,\dots , s -1\}$,  and
$$
 L= (A_0)^{\infty}= (A_0) \times  (A_0) \times  (A_0)\times\dots  
$$
be the space of one-sided sequences of  elements of $ A_0$.

An expansion of a real number $x \in [0,1]$ in the form
\begin{equation}
\label{eq: s-adic expansion}
x=\frac{\alpha_1}{s}+\frac{\alpha_2}{s^2}+\dots+\frac{\alpha_n}{s^n}+\dots ,
\end{equation}
where $\alpha_n \in A$, is called \emph{the s-adic expansion of $x$}.
By $x=\Delta^s _{\alpha_1\alpha_2...\alpha_n...}$ denote the s-adic expansion of $x$. The  notation $\Delta^s _{\alpha_1\alpha_2...\alpha_n...}$ is called \emph{the s-adic representation of $x$}.

Obviously, the  notation $x=\Delta^{-s} _{\alpha_1\alpha_2...\alpha_n...}$ is called \emph{the nega-s-adic representation of $x$}. Here
\begin{equation}
\label{eq: nega-s-adic expansion}
x=\Delta^{-s} _{\alpha_1\alpha_2...\alpha_n...}=-\frac{\alpha_1}{s}+\frac{\alpha_2}{s^2}-\frac{\alpha_3}{s^3}+\dots+\frac{(-1)^n\alpha_n}{s^n}+\dots ,
\end{equation}
where $\alpha_n \in A$.

 If $(k_n)$ is a certain fixed sequence of positive integers, then a series of the form 
$$
\frac{\alpha_{k_1}}{(-s)^{k_1}}+\frac{\alpha_{k_2}}{(-s)^{k_2}}+\dots+\frac{\alpha_{k_n}}{(-s)^{k_n}}+\dots , \alpha_{k_n} \in A,
$$
is {\itshape a nega-s-adic series.}

Suppose that 
$m_1=k_1$, $m_2=k_2-k_1$, $m_3=k_3-k_2$, \dots, $m_n=k_n-k_{n-1}$, \dots.
Then we obtain the following series
\begin{equation}
\label{ns1}
\sum^{\infty} _{n=1} {\frac{\alpha_{m_1+m_2+\dots+m_n}}{(-s)^{m_1+m_2+\dots+m_n}}}, 
\end{equation}
where $\alpha_{m_1+m_2+\dots+m_n} \in A.$

Numbers  $x \in \left[-\frac{s}{s+1},\frac{1}{s+1}\right]$ having a representation in form  \eqref{ns1} have the following nega-s-adic representation
$$
x=\sum^{\infty} _{n=1} {\frac{\alpha_{m_1+m_2+\dots+m_n}}{(-s)^{m_1+m_2+\dots+m_n}}}=\Delta^{-s} _{\underbrace{0\ldots 0}_{m_1-1}\alpha_{m_1}\underbrace{0\ldots 0}_{m_2-1}\alpha_{m_1+m_2}...\underbrace{0\ldots 0}_{m_n-1}\alpha_{m_1+m_2+\dots +m_n}...}.
$$

Let     $(d_n)$  be a fixed sequence of positive integers such that $d_n>1$ for all $n\in\mathbb N$, $(A_n)$ be a sequence of the sets  $A_n =\{0,1,2, \dots ,d_n-1\}$, and $L_n =A_1 \times A_2\times A_n\times \dots$.

A series of the form
\begin{equation}
\label{nc1}
-\frac{\varepsilon_1}{d_1}+\frac{\varepsilon_2}{d_1d_2}-\frac{\varepsilon_3}{d_1d_2d_3}+\dots+\frac{(-1)^n\varepsilon_n}{d_1d_2\cdots d_n}+\dots, 
\end{equation}
where $\varepsilon_n \in A_n$,  is called \emph{an alternating Cantor series}.

 In September 2013 (see the  presentation (in Ukrainian) and the working paper (in Ukrainian) that available at  https://www.researchgate.net/publication/303720347,  \\ https://www.researchgate.net/publication/316787375, respectively), the expansion of numbers by an alternating Cantor series was investigated as a numeral system,  and presented in the the report ``Representations of real numbers by alternating Cantor series"  at the fractal analysis  seminar of Institute of Mathematics of NAS of Ukraine and the National Pedagogical Dragomanov University. These results were published in  \cite{S. Serbenyuk alternating Cantor series 2013}.

An alternating Cantor series  that is a nega-s-adic series is called \emph{a nega-s-adic Cantor series}. That is 
\begin{equation}
\label{nsc1}
-\frac{\varepsilon_1}{s^{m_1}}+\frac{\varepsilon_2}{s^{m_1+m_2}}-\frac{\varepsilon_3}{s^{m_1+m_2+m_3}}+\dots +\frac{(-1)^n\varepsilon_n}{s^{m_1+m_2+\dots+m_n}}+\dots , \varepsilon_n \in A.
\end{equation}

It is easy to see that the following statement is true. 

\begin{lemma}[\cite{S. Serbenyuk 2013}] Nega-s-adic series  \eqref{ns1}  is an alternating Cantor series if and only if for any 
$n \in \mathbb N$ a sequence  $(m_n)$  is a sequence of odd positive integers and $\varepsilon_n=\alpha_n \in A$ as well. 
\end{lemma}

 A series of the form 
$$
-\frac{\alpha_{1}}{s^{k_1}}+\frac{\alpha_{2}}{s^{k_2}}-\frac{\alpha_{3}}{s^{k_3}}+\dots+\frac{(-1)^n\alpha_{n}}{s^{k_n}}+\dots , ~\alpha_{n} \in A.
$$
 is called \emph{a mixed s-adic series}.  Trivially, the last series is an alternating Cantor series.

We note that the case, when sequences $(\alpha_n)$ and $(m_n)$  are interdependent, is interesting, e.g. when $m_n=\alpha_n \in A_0$  for an arbitrary   $n \in \mathbb N$. In particular,  we shall describe properties of the set 
$$
S^{-}=\left\{x: x=\sum^{\infty} _{n=1} {\frac{(-1)^n\alpha_n}{s^{\alpha_1+\alpha_2+\dots+\alpha_n}}}, (\alpha_n) \in L, s>2 \right\}
$$
in the present article. Also, here the following set is considered:
$$
M_{(-D,s)}=\left\{x: x= \Delta^{-s} _{\underbrace{0\ldots 0}_{m_1-1}\alpha_{m_1}\underbrace{0\ldots 0}_{m_2-1}\alpha_{m_1+m_2}...\underbrace{0\ldots 0}_{m_n-1}\alpha_{m_1+m_2+\dots+m_n}...}\right\},
$$
where $s>1$ is a fixed positive integer, $\alpha_{m_1+m_2+\dots+m_n}\ne 0$ for all  $n \in \mathbb N$, and  $m_n\in \{3,5,7, \dots ,2i+1, \dots\}$.

\section{Fractal sets}

Let us consider the Cantor set. Any element of the Cantor set has only the digits $0$ and $2$
 in own ternary representation. This set  is  an  uncountable, perfect, and nowhere dense set of zero Lebesgue measure. Also, this is a self-similar fractal whose Hausdorff-Besicovitch dimension is equal to $\log_3 2$. 

One can formulate a general theorem about values of the Hausdorff-Besicovitch dimension of a set whose elements have 
restrictions on using combinations of digits in own s-adic  representation. 
\begin{theorem}[\cite{{Symon2}, {S. Serbenyuk 2017  fractals}}]%, {S. Serbenyuk 2013}]
\label{th: theorem6}
Let  $E$ be a set whose elements represented by a finite number of fixed combinations $\sigma_1, \sigma_2,\dots,\sigma_m$ of s-adic digits in the s-adic numeral system. Then the Hausdorff-Besicovitch dimension $\alpha_0$ of $E$ satisfies the following equation: 
$$
N(\sigma^1 _m)\left(\frac{1}{s}\right)^{\alpha_0}+N(\sigma^2 _m)\left(\frac{1}{s}\right)^{2\alpha_0}+\dots+N(\sigma^{k} _m)\left(\frac{1}{s}\right)^{k\alpha_0}=1,
$$
where $N(\sigma^k_m)$ is a number of k-digit combinations $\sigma^k_m$ from the set $\{\sigma_1, \sigma_2,\dots,\sigma_m\}$,
$k \in \mathbb N$, and $N(\sigma^1 _m)+N(\sigma^2 _m)+\dots+ N(\sigma^{k} _m)=m$.
\end{theorem}

This theorem is interesting since fractal properties of  many sets of special types follow from this theorem. For example,  the following set, whose elements have a functional restriction on using  digits in own the s-adic representation, was studied in \cite{Symon1}:
$$
S=\left\{ x: x=\sum^{\infty} _{n=1} {\frac{\alpha_n}{s^{\alpha_1+ \alpha_2+\dots+\alpha_n}}},  (\alpha_n) \in L\right\},
$$
where $s>2$ is a fixed positive integer. The last-mentioned  set is the set of all numbers whose s-adic representations contain only the following combinations of s-adic digits:
$$
1, 02, 003, \dots , \underbrace{0\ldots 0}_{i-1} {i}, \dots , \underbrace{0\ldots 0}_{s-2}[s-1].
$$ 
The Hausdorff-Besicovitch dimension $\alpha_0$ of the set $S$ satisfies the equation
$$
\left (\frac{1}{s}\right)^{\alpha_0}+\left (\frac{1}{s}\right)^{2\alpha_0}+\left (\frac{1}{s}\right)^{3\alpha_0}+\dots+\left (\frac{1}{s}\right)^{(s-1)\alpha_0}=1.
$$

Suppose  $s>2$ be a fixed positive integer number.

Consider a class  $\Upsilon_s$ of sets  $\mathbb  S_{(s,u)}$ represented  in the form 
\begin{equation*}
%\label{S(s,u)2}
\mathbb S_{(s,u)}= \left\{x: x=\frac{u}{s-1} +\sum^{\infty} _{n=1} {\frac{\alpha_n - u}{s^{\alpha_1+\dots+\alpha_n}}}, (\alpha_n) \in L, \alpha_n \ne u, \alpha_n \ne 0  \right\},
\end{equation*}
where $u=\overline{0,s-1}$, $u$ and $s$ are fixed for the set $\mathbb  S_{(s,u)}$. That is the class $\Upsilon_s$ contains the sets  $\mathbb  S_{(s,0)}, \mathbb  S_{(s,1)},\dots,\mathbb  S_{(s,s-1)}$. We say that   $\Upsilon$ is a class of sets such that contains the classes   $\Upsilon_3, \Upsilon_4,\dots ,\Upsilon_n,\dots$.

It is easy to see that the set  $\mathbb  S_{(s,u)}$ can be defined by the s-adic representation in the following form
\begin{equation*}
%\label{S(s,u)1}
\mathbb S_{(s,u)}=\left\{x: x= \Delta^{s}_{{\underbrace{u\ldots u}_{\alpha_1-1}} \alpha_1{\underbrace{u\ldots u}_{\alpha_2 -1}}\alpha_2 ...{\underbrace{u\ldots u}_{ \alpha_n -1}}\alpha_n...},  (\alpha_n) \in L, \alpha_n \ne u, \alpha_n \ne 0 \right\}, 
\end{equation*}

\begin{theorem}[\cite{{Symon2}, {S. Serbenyuk 2017  fractals}}]
For  an arbitrary $u \in A$ the set $\mathbb S_{(s,u)}$ is an uncountable,   perfect,   nowhere dense set of zero Lebesgue measure, and a self-similar fractal whose Hausdorff-Besicovitch dimension $\alpha_0 (\mathbb S_{(s,u)})$ satisfies the following equation 
$$
\sum _{p_i \ne u, p_i \in A_0} {\left(\frac{1}{s}\right)^{p_i \alpha_0}}=1.
$$
\end{theorem}

To prove the last statement, the auxiliary notion ``cylinder" is used. This notion is useful for study of local properties of considered sets (see the following lemma).

By $x_0=\Delta^{(s,u)} _{c_1...c_n...}$ denote the  equality 
$$
x_0=\frac{u}{s-1} +\sum^{\infty} _{k=1}{\frac{c_k - u}{s^{c_1+\dots+c_k}}}.
$$
That is
$$
 x_0 =\Delta^{(s,u)} _{c_1...c_n...} = \Delta^ {s}_{\underbrace{u...u}_{c_1 - 1}c_1 \underbrace{u...u}_{c_2 - 1}c_2...\underbrace{u...u}_{c_n - 1}c_n...}.
$$

\begin{definition} \emph{A cylinder   $ \Delta^{(s,u)} _{c_1...c_n}$ of rank $n$ with base $c_1c_2\ldots c_n$} is a set of the following
form
$$
\Delta^{(s,u)} _{c_1\ldots c_n}=\left\{x: x=\left(\sum^{n} _{k=1} {\frac{c_k -u}{s^{c_1+\dots+c_k}}}\right)+\frac{1}{s^{c_1+\dots+c_n}}{\left( \sum^{\infty} _{i=n+1} {\frac{\alpha_i - u}{s^{\alpha_{n+1}+\dots+\alpha_i}}}\right)}+\frac{u}{s-1}  \right\},
$$
 where $ c_1, c_2,\dots ,c_n $ are fixed s-adic digits, $c_n\ne 0$, $c_n\ne u$,  $\alpha_n \ne u,\alpha_n \ne 0$, and $ 2<s \in\mathbb N, n \in\mathbb N $.
\end{definition}

\begin{lemma}[\cite{Symon2, S. Serbenyuk 2017  fractals}] Cylinders $ \Delta^{(s,u)} _{c_1...c_n...} $ have the following properties:
\begin{enumerate}
\item
$$
\inf  \Delta^{(s,u)} _{c_1...c_n...}=\begin{cases}
\tau +\frac{1}{s^{c_1+...+c_n}}\left(\frac{s-1-u}{s^{s-1}-1}+\frac{u}{s-1}\right)&\text{if $u \in \{0,1\}$}\\
$$\\
\tau +\frac{1}{s^{c_1+\dots+c_n}}\frac{1}{s-1}&\text{if $ u \in \{2,3,\dots ,s-1\}$,}
\end{cases}
$$
$$
\sup  \Delta^{(s,u)} _{c_1...c_n...}=\begin{cases}
\tau + \frac{1}{s^{c_1+\dots+c_n}}\frac{1}{s-1}&\text{if $ u=0 $}
$$\\
$$\\
\tau +\frac{1}{s^{c_1+\dots+c_n}}\left(\frac{1}{s^{u+1}-1}+\frac{u}{s-1}\right)&\text{if $ u \in \{1,2,\dots ,s-2\} $}\\
$$\\
\tau +\frac{1}{s^{c_1+\dots+c_n}}\left(1-\frac{1}{s^{s-2}-1}\right)&\text{if $u=s-1$,}\\
\end{cases}
$$
where
$$
\tau= \sum^{n} _{k=1} {\frac{c_k-u}{s^{c_1+\dots+c_k}}}+\sum^{n} _{k=1} {\frac{u}{s^k}}.
$$
\item If $d(\cdot) $ is the diameter of a set, then
$$
d(\Delta^{(s,u)} _{c_1...c_n})=\frac{1}{s^{c_1+\dots+c_n}}d(\mathbb S_{(s,u)});
$$
\item
$$
\frac{d(\Delta^{(s,u)} _{c_1...c_nc_{n+1}})}{d(\Delta^{(s,u)} _{c_1...c_n})}=\frac{1}{s^{c_{n+1}}};
$$
\item 
$$
  \Delta^{(s,u)} _{c_1c_2...c_n} =\bigcup^{s-1} _{i=1} { \Delta^{(s,u)} _{c_1c_2...c_ni}}~~~\forall c_n \in A_0,~~~n \in \mathbb N,~ i \ne u.
$$
\item The following relationships hold: 
\begin{enumerate}
\item if $ u\in \{0,1\}$, then 
$$
\inf \Delta^{(s,u)} _{c_1...c_np}> \sup \Delta^{(s,u)} _{c_1...c_n[p+1]};
$$
\item if  $ u \in \{2,3,\dots ,s-3\}$, then 
$$
\begin{cases}
\sup \Delta^{(s,u)} _{c_1...c_np}< \inf \Delta^{(s,u)} _{c_1...c_n[p+1]}&\text{for all $p+1\le u$}\\
$$\\
\inf \Delta^{(s,u)} _{c_1...c_np}> \sup \Delta^{(s,u)} _{c_1...c_n[p+1]}&\text{for all $u<p$;}
\end{cases}
$$
\item if $ u  \in \{s-2,s-1\}$, then
$$
\sup \Delta^{(s,u)} _{c_1...c_np}< \inf \Delta^{(s,u)} _{c_1...c_n[p+1]}.
$$
\end{enumerate}
\end{enumerate}
\end{lemma}
The fifth property of the last lemma means the following: 
\begin{itemize}
\item for any positive integer $n$ cylinders $\Delta^{(s,u)} _{c_1...c_n}$ are right-to-left situated
in the case of the set $\mathbb S_{(s,0)}$ or $\mathbb S_{(s,1)}$;

\item let we have  the sets $\mathbb S_{(s,2)}$, $\mathbb S_{(s,3)}$, \dots, $\mathbb S_{(s,s-3)}$; then 
cylinders $\Delta^{(s,u)} _{c_1...c_n}$ ($u=\overline{2,s-3}$) are left-to-right situated for all $c_n \le 1$, $c_n \le 2$, \dots , $c_n \le s-4$, respectively, and cylinders $\Delta^{(s,u)} _{c_1...c_n}$ are right-to-left situated for all $c_n >2$, $c_n >3$, \dots , $c_n >s-3$,  respectively;

\item for all positive integers $n$ cylinders $\Delta^{(s,u)} _{c_1...c_n}$ are left-to-right situated
in the case of the set $\mathbb S_{(s,s-2)}$ or $\mathbb S_{(s,s-1)}$;

\item for any $\mathbb S_{(s,u)}$, $n\in\mathbb N$, and $c_n\ne s-1$ the following condition holds:
$$
\Delta^{(s,u)} _{c_1...c_{n-1}c_n}\cap\Delta^{(s,u)} _{c_1...c_{n-1}[c_n+1]}=\varnothing.
$$
\end{itemize} 
For proving the nowhere density of $\mathbb S_{(s,u)}$, the last property is used. 

Consider the set of all numbers whose s-adic representations contain only combinations of s-adic digits that is using  in the s-adic representations of elements of $\mathbb S_{(s,u)}$.

By $ \tilde S $ denote  the set of all numbers whose s-adic representations  contain only combinations of s-adic digits from the set
$$
\{1, 02, 003,\dots, {\underbrace{u\ldots u}_{c-1}}c, \dots, {\underbrace{(s-1)\ldots (s-1)}_{s-3}}(s-2)\},
$$
where $c \in A_0, u \in A$, $c \ne u$.

\begin{theorem}[\cite{Symon2, S. Serbenyuk 2017  fractals}]
The set $ \tilde S $ is:
\begin{itemize}
\item an uncountable, perfect, and nowhere dense set of zero Lebesgue measure;
\item a self-similar fractal, and its Hausdorff-Besicovitch dimension $\alpha_0$ satisfies the following equation 
$$
\left(\frac{1}{s}\right)^{\alpha_0}+(s-1)\left(\frac{1}{s}\right)^{2\alpha_0}+(s-1)\left(\frac{1}{s}\right)^{3\alpha_0}+\dots+(s-1)\left(\frac{1}{s}\right)^{(s-1)\alpha_0}=1.
$$
\end{itemize}
\end{theorem}

Let us prove the second item. 
The s-adic representation of an arbitrary element from $ \tilde S $ contains combinations of digits from the following  tuple:
$$
02, 003,\dots,\underbrace{0\ldots 00}_{s-2}(s-1);
$$
$$
1, 12, 113,\dots, \underbrace{1\ldots 11}_{s-2}(s-1);
$$
$$
223, 2224,\dots, \underbrace{2\ldots 22}_{s-2}(s-1);
$$
$$
\dots \dots \dots \dots \dots \dots \dots
$$
$$
u2, uu3,\dots,\underbrace{u\ldots uu}_{u-2}(u-1), \underbrace{u\ldots uu}_{u}(u+1), \dots, \underbrace{u\ldots uu}_{s-2}(s-1);
$$
$$
\dots \dots \dots \dots \dots \dots \dots
$$
$$
(s-1)2, (s-1)(s-1)3,\dots, \underbrace{(s-1)\ldots (s-1)(s-1)}_{s-3}(s-2).
$$
Here $s^2 -3s+3$ combinations of s-adic digits, i.e., the unique 1-digit combination  and $s-1$ k-digit combinations for all $k=\overline{2,s-1}$. Our statement follows from  Theorem \ref{th: theorem6}.

Let us consider some fractal sets defined in terms of the nega-s-adic representation, a nega-s-adic Cantor series, and a mixed s-adic series.

Let $s>2$ be a fixed positive integer.  

\begin{theorem}[\cite{S. Serbenyuk 2013}]
The sets 
$$
\mathbb S_{(-s,0)}=\left\{x: x=\sum^{\infty} _{n=1} {\frac{\alpha_n}{(-s)^{\alpha_1+\alpha_2+\dots+\alpha_n}}}, (\alpha_n) \in L\right\},
$$
$$
S^{-}=\left\{x: x=\sum^{\infty} _{n=1} {\frac{(-1)^n\alpha_n}{s^{\alpha_1+\alpha_2+\dots+\alpha_n}}}, (\alpha_n) \in L \right\},
$$
 are:
\begin{itemize}
\item	 uncountable,   perfect,   nowhere dense sets of zero Lebesgue measure;
\item	self-similar fractals whose Hausdorff-Besicovitch dimension $\alpha_0$ satisfies the following equation 
 $$
\sum^{s-1} _{i=1} {\left(\frac{1}{s}\right)^{i\alpha_0}}=1.
$$
\end{itemize}
\end{theorem}
\begin{proof} Let us prove that the sets $ \mathbb S_{(-s,0)}$ and $S^{-}$ are uncountable.

Let us prove that the sets $\mathbb S_{(-s,0)}$ and $C[-s, A_0]$ are equivalent. That is, let us consider the mapping   
$$
x=\sum^{\infty} _{n=1} {\frac{\alpha_n \cdot (-1)^{\alpha_1+\alpha_2+...+\alpha_n} }{s^{\alpha_1+\alpha_2+...+\alpha_n}}} \stackrel{f}{\longrightarrow} \sum^{\infty} _{n=1} {\frac{\alpha_n}{(-s)^n}}=f(x)=y
$$
or (in other words)
$$
 x=\Delta^{-s} _{\underbrace{0...0}_{\mbox{$\alpha_1-1$ }}\alpha_1 \underbrace{0...0}_{\mbox{$\alpha_2-1$ }}\alpha_2...\underbrace{0...0}_{\mbox{$\alpha_n-1$ }}\alpha_n...}  \stackrel{f}{\longrightarrow}  \Delta^{-s} _{\alpha_1\alpha_2...\alpha_n...}=f(x)=y.
$$

Suppose $x_1$ and  $x_2$ from $\mathbb S_{(-s,0)}$ are such that $x_1 \ne x_2$ and
$$
x_1=\Delta^{-s} _{\underbrace{0...0}_{\mbox{$\alpha_1-1$ }}\alpha_1 \underbrace{0...0}_{\mbox{$\alpha_2-1$ }}\alpha_2...\underbrace{0...0}_{\mbox{$\alpha_n-1$ }}\alpha_n...}, ~~~x_2=\Delta^{-s} _{\underbrace{0...0}_{\mbox{$\beta_1-1$ }}\beta_1 \underbrace{0...0}_{\mbox{$\beta_2-1$ }}\beta_2...\underbrace{0...0}_{\mbox{$\beta_n-1$ }}\beta_n...}.
$$

If $f(x_1)=f(x_2)$ is nega-s-adic irrational (i.e., this number has the unique representation).   Hence $\alpha_n=\beta_n$ holds for all $n \in \mathbb N$. That is, $x_1=x_2$, It contradicts to the condition. 

Assume that $f(x_1)=f(x_2)$ is nega-s-adic rational, but this is not possible because any number from  $C[-s, A_0]$ does not have two expansions. 

So, $f$ is a bijection. Since $C[-s, A_0]$ is a uncountable set, we
see that  $\mathbb S_{(-s,0)}$ is a uncountable set. Proofs for $S^-$ are similar.

Statements of this theorem follows from properties of the following notions of cylinders. Proofs are similar with proofs of  Theorems 1 and 3 in \cite{{S. Serbenyuk 2017  fractals}} (arXiv:1703.05262). 
\end{proof}
\begin{definition}
\emph{A cylinder  $\Delta^{(-s, 0)} _{c_1c_2...c_n}$ of rank $n$ with base $c_1c_2\ldots c_n$}  is a set 
formed by all numbers of the set $\mathbb S_{(-s,0)}$ with nega-s-representations in
which the first $n$ non-zero digits coincide with $c_1, c_2, \dots , c_n$ respectively.
\end{definition}

\begin{definition}
\emph{ A cylinder $\Delta^{-} _{c_1c_2...c_n}$ of rank $n$ with base $c_1c_2\ldots c_n$}  is a subset of $S^{-}$ with elemets for which the following condition holds:
$$
\alpha_1=c_1, \alpha_2=c_2,\dots , \alpha_n=c_n, 
$$
where  $c_1,c_2,\dots ,c_n$ is  an ordered tuple of numbers. 
\end{definition}

\begin{lemma}[\cite{S. Serbenyuk 2013}]
Cylinders $\Delta^{(-s, 0)} _{c_1c_2...c_n}$ have the following properties:
\begin{enumerate}
\item
$$
\inf \Delta^{(-s, 0)} _{c_1c_2...c_n}=\begin{cases}
g^{(-s)} _n+ \frac{\inf \mathbb S_{(-s,0)}}{(-s)^{c_1+c_2+\dots+c_n}}&\text{if  $c_1+\dots+c_n$ is even}\\
g^{(-s)} _n+ \frac{\sup \mathbb S_{(-s,0)}}{(-s)^{c_1+c_2+\dots+c_n}}&\text{if  $c_1+\dots+c_n$ is odd,}
\end{cases}
$$
$$
\sup \Delta^{(-s, 0)} _{c_1c_2...c_n}=\begin{cases}
g^{(-s)} _n+ \frac{\sup \mathbb S_{(-s,0)}}{(-s)^{c_1+c_2+\dots+c_n}}&\text{if  $c_1+\dots+c_n$ is even}\\
g^{(-s)} _n+ \frac{\inf \mathbb S_{(-s,0)}}{(-s)^{c_1+c_2+\dots+c_n}}&\text{if  $c_1+\dots+c_n$ is odd,}
\end{cases}
$$
where
$$
g^{(-s)} _n=\sum^n _{i=1} {\frac{c_i(-1)^i}{s^{c_1+c_2+\dots+c_i}}}.
$$
\item Suppose $d(\cdot)$ is the diameter of a set. Then 
$$
d\left(\Delta^{(-s, 0)} _{c_1c_2...c_n}\right)=\frac{d\left(\mathbb S_{(-s,0)}\right)}{s^{c_1+c_2+\dots+c_n}}.
$$
\item The main metric relationship is following:
$$
\frac{d\left(\Delta^{(-s, 0)} _{c_1c_2...c_nc_{n+1}}\right)}{d\left(\Delta^{(-s, 0)} _{c_1c_2...c_n}\right)}=\frac{1}{s^{c_{n+1}}}.
$$
\item For any $n\in\mathbb N$ the following condition holds:
$$
\Delta^{(-s, 0)} _{c_1c_2...c_n}=\bigcup^{s-1} _{i=1} {\Delta^{(-s, 0)} _{c_1c_2...c_ni}}.
$$
\item For cylinders  $\Delta^{(-s, 0)} _{c_1c_2...c_nc_{n+1}}$ of rank $(n+1)$ with base $c_1c_2\dots c_nc_{n+1}$ the following relationships hold:
$$
\inf \Delta^{(-s, 0)} _{c_1c_2...c_np} > \sup \Delta^{(-s, 0)} _{c_1c_2...c_n [p+1]} ~ \mbox{whenever} ~ c_1+c_2+\dots+c_n+p \mbox{ is even},
$$
$$
~~\inf \Delta^{(-s, 0)} _{c_1c_2...c_n [p+1]} > \sup \Delta^{(-s, 0)} _{c_1c_2...c_np}~ \mbox{whenever} ~ c_1+c_2+\dots+c_n+p \mbox{ is odd}.
$$
\item Let $T^{(-s,0)} _{c_1c_2...c_np} $ be an interval of the form
$$
T^{(-s,0)} _{c_1c_2...c_np}=\begin{cases}
\left(\sup \Delta^{(-s, 0)} _{c_1c_2...c_n [p+1]}, \inf \Delta^{(-s, 0)} _{c_1c_2...c_n p}\right)&\text{if $c_1+\dots+c_n+p$ is even}\\
\left(\sup \Delta^{(-s, 0)} _{c_1c_2...c_np }, \inf \Delta^{(-s, 0)} _{c_1c_2...c_n [p+1]}\right)&\text{if $c_1+\dots+c_n+p$ is odd,}
\end{cases}
$$
where $1 \le p< s-1$ is a positive integer. Then
$$
T^{(-s,0)} _{c_1c_2...c_np} \cap \mathbb S_{(-s,0)}= \varnothing.
$$
\item For any $p \in \{1, 2, \dots , s-2\}$ the following condition holds: 
$$
\Delta^{(-s, 0)} _{c_1c_2...c_np } \cap \Delta^{(-s, 0)} _{c_1c_2...c_n [p+1]}=\varnothing.
$$
\item If $x_0 \in \mathbb S_{(-s,0)}$, then
$$
x_0=\bigcap^{\infty} _{n=1} {\Delta^{(-s, 0)} _{c_1c_2...c_n}}.
$$
\end{enumerate}
\end{lemma}
\begin{proof} {\itshape Properties 1 -- 4} follow from the definitions of 
$\Delta^{(-s, 0)} _{c_1c_2...c_n}$ and $\mathbb S_{(-s,0)}$.

Let us prove {\itshape Property 5}. Suppose $\Delta^{(-s, 0)} _{c_1c_2...c_np},  \Delta^{(-s, 0)} _{c_1c_2...c_n (p+1)}$ are cylinders, where $1 \le p<s-1$, and
$$
g^{(-s)} _n=\sum^{n} _{i=1} {\frac{c_i}{(-s)^{c_1+c_2+...+c_i}}};~~~\varpi_n=c_1+c_2+...+c_n.
$$

From the definition of  $\Delta^{(-s, 0)} _{c_1c_2...c_n}$, it follows that  
$$
 \Delta^{(-s, 0)} _{c_1c_2...c_np} \subset \begin{cases}
\left[g^{(-s)} _n+ \frac{p}{(-s)^{\varpi_n+p}}+\frac{-{(s^2+1)}}{s(s^2-1)(-s)^{\varpi_n+p}};g^{(-s)} _n+ \frac{p}{(-s)^{\varpi_n+p}}+\frac{2}{(s^2-1)(-s)^{\varpi_n+p}}\right],\\
\left[g^{(-s)} _n+ \frac{p}{(-s)^{\varpi_n+p}}+\frac{2}{(s^2-1)(-s)^{\varpi_n+p}};g^{(-s)} _n+ \frac{p}{(-s)^{\varpi_n+p}}+\frac{-{(s^2+1)}}{s(s^2-1)(-s)^{\varpi_n+p}}\right],
\end{cases}
$$
where $\varpi_n+p$ is even for the first case, and $\varpi_n+p$ is odd for the second case.

By analogy, we obtain
$$
\Delta^{(-s, 0)} _{c_1c_2...c_n(p+1)} \subset \begin{cases}
\left[g^{(-s)} _n+ \frac{p+1}{(-s)^{\varpi_n+p+1}}+\frac{-{(s^2+1)}}{s(s^2-1)(-s)^{\varpi_n+p+1}};g^{(-s)} _n+ \frac{p+1}{(-s)^{\varpi_n+p+1}}+\frac{2}{(s^2-1)(-s)^{\varpi_n+p+1}}\right],\\
\left[g^{(-s)} _n+ \frac{p+1}{(-s)^{\varpi_n+p+1}}+\frac{2}{(s^2-1)(-s)^{\varpi_n+p+1}};g^{(-s)} _n+ \frac{p+1}{(-s)^{\varpi_n+p+1}}+\frac{-{(s^2+1)}}{s(s^2-1)(-s)^{\varpi_n+p+1}}\right],
\end{cases}
$$
where  $(\varpi_n+p+1)$  is even for the first case and is odd  for the second case.

Let us prove the mentioned inequalities. 

Let $\varpi_n+p=c_1+c_2+...+c_n+p$ be an even number. Then  
$$
\inf \Delta^{(-s, 0)} _{c_1c_2...c_np} - \sup \Delta^{(-s, 0)} _{c_1c_2...c_n (p+1)}=g^{(-s)} _n+ \frac{p}{(-s)^{\varpi_n+p}}+\frac{-{(s^2+1)}}{s(s^2-1)(-s)^{\varpi_n+p}}-
$$
$$
- g^{(-s)} _n- \frac{p+1}{(-s)^{\varpi_n+p+1}}-\frac{-{(s^2+1)}}{s(s^2-1)(-s)^{\varpi_n+p+1}}=\frac{1}{s^{\varpi_n+p}}\left(ps+p+1-\frac{s^3+s^2+s+1}{s(s^2-1)}\right)>0,
$$
because
$$
\frac{s^3+s^2+s+1}{s(s^2-1)}=1+\frac{(s+1)^2}{s(s^2-1)}=1+\frac{s+1}{s(s-1)}<2.
$$

Let $\varpi_n+p=c_1+c_2+...+c_n+p$ be an odd number. Then
$$
\inf \Delta^{(-s, 0)} _{c_1c_2...c_n (p+1)} - \sup \Delta^{(-s, 0)} _{c_1c_2...c_np}=g^{(-s)} _n+ \frac{p+1}{(-s)^{\varpi_n+p+1}}+\frac{-{(s^2+1)}}{s(s^2-1)(-s)^{\varpi_n+p+1}}
$$
$$
- g^{(-s)} _n- \frac{p}{(-s)^{\varpi_n+p}}-\frac{-{(s^2+1)}}{s(s^2-1)(-s)^{\varpi_n+p}}=\frac{1}{s^{\varpi_n+p+1}}\left(ps+p+1-\frac{s^3+s^2+s+1}{s(s^2-1)}\right)>0.
$$

To prove  {\itshape Property 6}, it suffices to prove the following inequalities:
\begin{itemize}
\item under the condition that  $c_1+c_2+...+c_n+p$ is an even number
$$
\left\{
\begin{array}{rcl}
\sup \Delta^{(-s,0)} _{c_1c_1...c_n(p+1)c_{n+2}}-\sup \Delta^{(-s,0)} _{c_1c_1...c_n(p+1)}&<&0,\\
\\
\inf \Delta^{(-s,0)} _{c_1c_1...c_npc_{n+2}}-\inf \Delta^{(-s,0)} _{c_1c_1...c_np} & > &0.\\
\end{array}
\right.
$$
\item under the condition that $c_1+c_2+...+c_n+p$ is an odd number
$$
\left\{
\begin{array}{rcl}
\sup \Delta^{(-s,0)} _{c_1c_1...c_npc_{n+2}}-\sup \Delta^{(-s,0)} _{c_1c_1...c_np}&<&0,\\
\\
\inf \Delta^{(-s,0)} _{c_1c_1...c_n(p+1)c_{n+2}}-\inf \Delta^{(-s,0)} _{c_1c_1...c_n(p+1)} & > &0.\\
\end{array}
\right.
$$
\end{itemize}

Suppose
$$
l_0 (c_1,c_2,...,c_n,p)=\begin{cases}
-\frac{s^2+1}{s(s^2-1)} &\text{whenever $c_1+c_2+...+c_n+p$ is even}\\
\frac{2}{s^2-1} &\text{whenever $c_1+c_2+...+c_n+p$ is odd,}
\end{cases}
$$
$$
l (c_1,c_2,...,c_n,p)=\begin{cases}
\frac{2}{s^2-1} &\text{whenever $c_1+c_2+...+c_n+p$ is even}\\
-\frac{s^2+1}{s(s^2-1)} &\text{whenever $c_1+c_2+...+c_n+p$ is odd.}
\end{cases}
$$

Let $c_1+c_2+...+c_n+p$ be an even number. Then 
$$
\sup \Delta^{(-s,0)} _{c_1c_1...c_n(p+1)c_{n+2}}-\sup \Delta^{(-s,0)} _{c_1c_1...c_n(p+1)}=g^{(-s)} _n+\frac{p+1}{(-s)^{c_1+...+c_n+p+1}}+\frac{c_{n+2}}{(-s)^{c_1+...+c_n+p+1+c_{n+2}}}+
$$
$$
+\frac{l (c_1,c_2,...,c_n,p+1, c_{n+2})}{(-s)^{c_1+...+c_n+p+1+c_{n+2}}}-g^{(-s)} _n-\frac{p+1}{(-s)^{c_1+...+c_n+p+1}}-\frac{l (c_1,c_2,...,c_n,p+1)}{(-s)^{c_1+...+c_n+p+1}}=
$$
$$
=-\frac{1}{s^{c_1+...+c_n+p+1}}\left(\frac{c_{n+2}}{(-s)^c_{n+2}}+\frac{l (c_1,c_2,...,c_n,p+1, c_{n+2})}{(-s)^c_{n+2}}+\frac{s^2+1}{s(s^2-1)}\right)=
$$
$$
=\begin{cases}
-\frac{1}{s^{c_1+...+c_n+p+1}}\left(\frac{c_{n+2}}{s^{c_n+2}}-\frac{s^2+1}{s(s^2-1)s^{c_{n+2}}}+\frac{s^2+1}{s(s^2-1)}\right)<0 &\text{if $c_{n+2}$ is even}\\
\\
-\frac{1}{s^{c_1+...+c_n+p+1}}\left(-\frac{c_{n+2}}{s^{c_n+2}}-\frac{2}{(s^2-1)s^{c_{n+2}}}+\frac{s^2+1}{s(s^2-1)}\right)<0 &\text{if  $c_{n+2}$ is odd,}
\end{cases}
$$
because
$$
-\frac{c_{n+2}}{s^{c_n+2}}-\frac{2}{(s^2-1)s^{c_{n+2}}}+\frac{s^2+1}{s(s^2-1)}=\frac{(s^2+1)s^{c_{n+2}}+sc_{n+2}-(s^2c_{n+2}+2)s}{(s^2-1)s^{1+c_{n+2}}}\ge 0.
$$

By analogy, we have
$$
\inf \Delta^{(-s,0)} _{c_1c_1...c_npc_{n+2}}-\inf \Delta^{(-s,0)} _{c_1c_1...c_np}=\frac{c_{n+2}}{(-s)^{c_1+...+c_n+p+c_{n+2}}}+\frac{l_0 (c_1,c_2,...,c_n,p, c_{n+2})}{(-s)^{c_1+...+c_n+p+c_{n+2}}}-
$$
$$
-\frac{l_0 (c_1,c_2,...,c_n,p)}{(-s)^{c_1+...+c_n+p}}=\begin{cases}
\frac{1}{s^{c_1+...+c_n+p}}\left(\frac{c_{n+2}}{s^{c_n+2}}-\frac{s^2+1}{s(s^2-1)s^{c_{n+2}}}+\frac{s^2+1}{s(s^2-1)}\right)>0 &\text{ if $c_{n+2}$ is even}\\
\\
\frac{1}{s^{c_1+...+c_n+p}}\left(-\frac{c_{n+2}}{s^{c_n+2}}-\frac{2}{(s^2-1)s^{c_{n+2}}}+\frac{s^2+1}{s(s^2-1)}\right)>0 &\text{if $c_{n+2}$ is odd.}
\end{cases}
$$

Let $c_1+c_2+...+c_n+p$ be an odd number.  Then
$$
\sup \Delta^{(-s,0)} _{c_1c_1...c_npc_{n+2}}-\sup \Delta^{(-s,0)} _{c_1c_1...c_np}=\frac{c_{n+2}}{(-s)^{c_1+...+c_n+p+c_{n+2}}}+\frac{l (c_1,c_2,...,c_n,p, c_{n+2})}{(-s)^{c_1+...+c_n+p+c_{n+2}}}-
$$
$$
-\frac{l (c_1,c_2,...,c_n,p)}{(-s)^{c_1+...+c_n+p}}=\begin{cases}
-\frac{1}{s^{c_1+...+c_n+p}}\left(\frac{c_{n+2}}{s^{c_n+2}}-\frac{s^2+1}{s(s^2-1)s^{c_{n+2}}}+\frac{s^2+1}{s(s^2-1)}\right)<0 &\text{if $c_{n+2}$ is even}\\
\\
-\frac{1}{s^{c_1+...+c_n+p}}\left(-\frac{c_{n+2}}{s^{c_n+2}}-\frac{2}{(s^2-1)s^{c_{n+2}}}+\frac{s^2+1}{s(s^2-1)}\right)<0 &\text{if $c_{n+2}$ is odd.}
\end{cases}
$$
Also,
$$
\inf \Delta^{(-s,0)} _{c_1c_1...c_n(p+1)c_{n+2}}-\inf \Delta^{(-s,0)} _{c_1c_1...c_n(p+1)}=\frac{c_{n+2}}{(-s)^{c_1+...+c_n+p+1+c_{n+2}}}+\frac{l_0(c_1,c_2,...,c_n,p+1, c_{n+2})}{(-s)^{c_1+...+c_n+p+1+c_{n+2}}}-
$$
$$
-\frac{l_0(c_1,c_2,...,c_n,p+1)}{(-s)^{c_1+...+c_n+p+1}}=\begin{cases}
\frac{1}{s^{c_1+...+c_n+p}}\left(\frac{c_{n+2}}{s^{c_n+2}}-\frac{s^2+1}{s(s^2-1)s^{c_{n+2}}}+\frac{s^2+1}{s(s^2-1)}\right)>0 &\text{if $c_{n+2}$ is even}\\
\\
\frac{1}{s^{c_1+...+c_n+p}}\left(-\frac{c_{n+2}}{s^{c_n+2}}-\frac{2}{(s^2-1)s^{c_{n+2}}}+\frac{s^2+1}{s(s^2-1)}\right)>0 &\text{if $c_{n+2}$ is odd.}
\end{cases}
$$

{\itshape Property 7} follows from   {Property 6}.

{\itshape Property 8}. From properties of cylinders of $\mathbb S_{(-s,0)}$, it follows the following:  if $x_0 \in \mathbb S_{(-s,0)}$, then
$$
x_0 \in \Delta^{(-s,0)} _{\alpha_1} \cap \Delta^{(-s,0)} _{\alpha_1\alpha_2} \cap ... \cap \Delta^{(-s,0)} _{\alpha_1\alpha_2...\alpha_n} \cap ...,
$$
where $x_0=\Delta^{-s} _{\underbrace{0...0}_{\mbox{$\alpha_1-1$ }}\alpha_1 \underbrace{0...0}_{\mbox{$\alpha_2-1$ }}\alpha_2...\underbrace{0...0}_{\mbox{$\alpha_n-1$ }}\alpha_n...} $.
Also,
$$
x_0 \in \left[\inf \Delta^{(-s,0)} _{\alpha_1}; \sup \Delta^{(-s,0)} _{\alpha_1} \right] \cap \left[\inf \Delta^{(-s,0)} _{\alpha_1\alpha_2}; \sup \Delta^{(-s,0)} _{\alpha_1\alpha_2} \right] \cap ...\cap \left[\inf \Delta^{(-s,0)} _{\alpha_1\alpha_2...\alpha_n}; \sup \Delta^{(-s,0)} _{\alpha_1\alpha_2...\alpha_n} \right] \cap ....
$$

So,  $x_0$ belongs to the following system of closed intervals:
$$
\left[\inf \Delta^{(-s,0)} _{\alpha_1}; \sup \Delta^{(-s,0)} _{\alpha_1} \right] \supset \left[\inf \Delta^{(-s,0)} _{\alpha_1\alpha_2}; \sup \Delta^{(-s,0)} _{\alpha_1\alpha_2} \right] \supset ... \supset \left[\inf \Delta^{(-s,0)} _{\alpha_1\alpha_2...\alpha_n}; \sup \Delta^{(-s,0)} _{\alpha_1\alpha_2...\alpha_n} \right]   \supset ....
$$
Therefore,
$$
x_0=\bigcap^{\infty} _{n=1} {\Delta^{(-s, 0)} _{c_1c_2...c_n}}.
$$
\end{proof}

\begin{lemma}[\cite{S. Serbenyuk 2013}]
Cylinders $\Delta^{-} _{c_1c_2...c_n} $  have the following properties: 
\begin{enumerate}
\item
$$
\Delta^{-} _{c_1c_2...c_n} \subset \begin{cases}
\left[\sigma_{2k}+\frac{\inf S^{-}}{s^{c_1+c_2+\dots+c_{2k}}},\sigma_{2k}+\frac{\sup S^{-}}{s^{c_1+c_2+\dots+c_{2k}}}\right]&\text{if $n=2k$ }\\

\left[\sigma_{2k+1}-\frac{\sup S^{-}}{s^{c_1+c_2+\dots+c_{2k+1}}},\sigma_{2k+1}-\frac{\inf S^{-}}{s^{c_1+c_2+\dots+c_{2k+1}}}\right]&\text{if $n=2k+1$,}
\end{cases}
$$
where $k \in \mathbb N$,
$$
\sigma_n=\sum^{n} _{i=1} {\frac{c_i}{s^{c_1+c_2+\dots+c_i}}}, ~\inf S^{-}=\frac{-s^{s-1}+s-1}{s^s-1}, ~\sup S^{-}=\frac{-s^2+s+1}{s^s-1}.
$$
\item
$$
d(\Delta^{-} _{c_1c_2...c_n})=\frac{s^{s-1}-s^2+2}{(s^s-1)s^{c_1+c_2+...+c_n}}.
$$
\item
$$
\frac{\Delta^{-} _{c_1c_2...c_nc_{n+1}}}{\Delta^{-} _{c_1c_2...c_n}}=\frac{1}{s^{c_{n+1}}}.
$$
\item
$$
\Delta^{-} _{c_1c_2...c_nc_{n+1}} \subset \Delta^{-} _{c_1c_2...c_n}~ \forall c_n \in A_0, ~n \in \mathbb N.
$$
\item
Cylinders $\Delta^{-} _{c_1c_2...c_{n-1}1}, \Delta^{-} _{c_1c_2...c_{n-1}2},\dots ,\Delta^{-} _{c_1c_2...c_{n-1}[s-1]}$ are:
\begin{itemize}
\item   right-to-left situated whenever  $n$ is even, i.e., 
$$
\forall k \in \mathbb N: \sup \Delta^{-} _{c_1c_2...c_{2k-1}[c_{2k}+1]}<\inf\Delta^{-} _{c_1c_2...c_{2k-1}c_{2k}};
$$
\item left-to-right situated whenever  $n$ is odd, i.e.,  
$$
\forall k \in \mathbb N: \sup \Delta^{-} _{c_1c_2...c_{2k}c_{2k+1}}<\inf\Delta^{-} _{c_1c_2...c_{2k}[c_{2k+1}+1]}.
$$
\end{itemize}
\end{enumerate}
\end{lemma}
\begin{proof}
{\itshape The first, the second}, and  {\itshape the  third} properties follow from  the definition of a cylinder $\Delta^{-} _{c_1c_2...c_n}$. 

Let us prove  {\itshape the fourth} property. 

1. Suppose $n=2k, k \in \mathbb N$. Then the equality $\inf \Delta^{-} _{c_1c_2...c_nc_{n+1}} \ge \inf\Delta^{-} _{c_1c_2...c_n}$ can be written in the form  
$$
\sum^{2k} _{m=1} {\frac{(-1)^mc_m}{s^{c_1+c_2+...+c_m}}}-\frac{c_{2k+1}}{s^{c_1+c_2+...+c_{2k+1}}}-\frac{\sup S^{-}}{s^{c_1+c_2+...+c_{2k+1}}}\ge\sum^{2k} _{m=1} {\frac{(-1)^mc_m}{s^{c_1+c_2+...+c_m}}}+\frac{\inf S^{-}}{s^{c_1+c_2+...+c_{2k}}}
$$
or
$$
-c_{2k+1}-\sup S^{-}\ge s^{c_{2k+1}}\inf S^{-},
$$
$$
\frac{(s^2-s-1)+s^{c_{2k+1}}(s^{s-1}-s+1)-c_{2k+1}(s^s-1)}{s^s-1}\ge0.
$$

It is easy to see that the last inequality is an equality under the condition  $c_{2k+1}=~1$.

In addition, for an even number $n$, let us consider the inequality 
$\sup \Delta^{-} _{c_1c_2...c_nc_{n+1}} \le \sup\Delta^{-} _{c_1c_2...c_n}$.
We have
$$
\sum^{2k} _{m=1} {\frac{(-1)^mc_m}{s^{c_1+c_2+...+c_m}}}-\frac{c_{2k+1}}{s^{c_1+c_2+...+c_{2k+1}}}-\frac{\inf S^{-}}{s^{c_1+c_2+...+c_{2k+1}}}\le\sum^{2k} _{m=1} {\frac{(-1)^mc_m}{s^{c_1+c_2+...+c_m}}}+\frac{\sup S^{-}}{s^{c_1+c_2+...+c_{2k}}},
$$
$$
-c_{2k+1}-\inf S^{-}\le s^{c_{2k+1}}\sup S^{-}
$$
or
$$
(1+c_{2k+1}+s^{2+c_{2k+1}}+s^{s-1})-s-s^{1+c_{2k+1}}-c_{2k+1}s^s-s^{c_{2k+1}}\le 0.
$$
The last inequality is an equality when $c_{2k+1}=s-1$ holds.

2. Suppose $n=2k+1, k \in \mathbb N$. Then $\inf \Delta^{-} _{c_1c_2...c_nc_{n+1}} \ge \inf\Delta^{-} _{c_1c_2...c_n}$ and
$$
\sum^{2k+1} _{m=1} {\frac{(-1)^mc_m}{s^{c_1+c_2+...+c_m}}}+\frac{c_{2k+2}}{s^{c_1+c_2+...+c_{2k+2}}}+\frac{\inf S^{-}}{s^{c_1+c_2+...+c_{2k+2}}}\ge\sum^{2k+1} _{m=1} {\frac{(-1)^mc_m}{s^{c_1+c_2+...+c_m}}}-\frac{\sup S^{-}}{s^{c_1+c_2+...+c_{2k+1}}}
$$
are equivalent. Hence 
$$
(s-1-c_{2k+2})-s^{c_{2k+2}}(s^2-s-1)+s^{s-1}(sc_{2k+2}-1)\ge 0.
$$
 If $c_{2k+2}=s-1$, then the last inequality is an equlity.  

By analogy, for $\sup \Delta^{-} _{c_1c_2...c_nc_{n+1}} \le \sup\Delta^{-} _{c_1c_2...c_n}$, we get
$$
\sum^{2k+1} _{m=1} {\frac{(-1)^mc_m}{s^{c_1+c_2+...+c_m}}}+\frac{c_{2k+2}}{s^{c_1+c_2+...+c_{2k+2}}}+\frac{\sup S^{-}}{s^{c_1+c_2+...+c_{2k+2}}}\le\sum^{2k+1} _{m=1} {\frac{(-1)^mc_m}{s^{c_1+c_2+...+c_m}}}-\frac{\inf S^{-}}{s^{c_1+c_2+...+c_{2k+1}}},
$$
$$
\sup S^{-}+c_{2k+2}(s^s-1)\le -s^{c_{2k+2}}\inf S^{-},
$$
and
$$
(s-s^2)+(1-c_{2k+2})+(s-1)s^{c_{2k+2}}+s^{s-1}(sc_{2k+2}-s^{c_{2k+2}})\le 0.
$$
It is true for all  values of $c_{2k+2}$ and  $s>2$, and is an equality when  $c_{2k+2}=1$.

Let us prove \emph{Property 5}.
\begin{itemize}
\item 
$$
\forall k \in \mathbb N: ~ \sup \Delta^{-} _{c_1c_2...c_{2k-1}[c_{2k}+1]}- \inf \Delta^{-} _{c_1c_2...c_{2k-1}c_{2k}}=
$$
$$
=\frac{c_{2k}+1}{s^{c_1+c_2+...+c_{2k}+1}}+\frac{\sup S^{-}}{{s^{c_1+c_2+...+c_{2k}+1}}}-\frac{c_{2k}}{{s^{c_1+c_2+...+c_{2k}}}}-\frac{\inf S^{-}}{{s^{c_1+c_2+...+c_{2k}}}}=
$$
$$
=\frac{1}{{s^{c_1+c_2+...+c_{2k}}}}\left(\frac{1-s}{s}c_{2k}+2\frac{s^s-s^2+s}{s(s^s-1)}\right)=
$$
$$
=\frac{s^s(2+c_{2k}-sc_{2k})+s(c_{2k}+2-2s)-c_{2k}}{s(s^s-1)s^{c_1+c_2+...+c_{2k}}}<0.
$$
\item
$$
\forall k \in \mathbb N: ~ \sup \Delta^{-} _{c_1c_2...c_{2k}c_{2k+1}}- \inf \Delta^{-} _{c_1c_2...c_{2k}[c_{2k+1}+1]}=
$$
$$
=\frac{c_{2k+1}}{s^{c_1+c_2+...+c_{2k}+c_{2k+1}}}-\frac{\inf S^{-}}{s^{c_1+c_2+...+c_{2k}+c_{2k+1}}}-\frac{1+c_{2k+1}}{s^{1+c_1+c_2+...+c_{2k}+c_{2k+1}}}+\frac{\sup S^{-}}{s^{1+c_1+c_2+...+c_{2k}+c_{2k+1}}}=
$$
$$
=\frac{1}{s^{1+c_1+c_2+...+c_{2k}+c_{2k+1}}}\left(\frac{s^s-2s^2+2s+1}{s^s-1}-\frac{s^{s+1}c_{2k+1}-sc_{2k+1}}{s^s-1}+\frac{s^sc_{2k+1}-c_{2k+1}}{s^s-1}+\frac{s^s-1}{s^s-1}\right)=
$$
$$
=\frac{s^s(2+c_{2k+1}-sc_{2k+1})+s(2-2s-c_{2k+1})-c_{2k+1}}{(s^s-1)s^{1+c_1+c_2+...+c_{2k}+c_{2k+1}}}<0.
$$
\end{itemize}
\end{proof}

It follows from the last-mentioned lemma that the following statements are true.
\begin{corollary}
 For all $c_n \in \{1,2,\dots ,s-2\}$ the condition 
$$
\Delta^{-} _{c_1c_2...c_{n-1}c_n} \cap \Delta^{-} _{c_1c_2...c_{n-1}[c_n+1]}=\varnothing
$$
 holds.
\end{corollary}

\begin{corollary}
Intervals of the form 
$$
\left(\sup \Delta^{-} _{c_1c_2...c_{2k-1}}, \inf \Delta^{-} _{c_1c_2...c_{2k-2}[c_{2k-1}+1]}\right)~\text{and}~\left(\sup \Delta^{-} _{c_1c_2...c_{2k-1}[c_{2k}+1]}, \inf \Delta^{-} _{c_1c_2...c_{2k-1}c_{2k}}\right), 
$$
where $k \in \mathbb N$, have the empty intersection with the set $S^{-}$.
\end{corollary}

\begin{corollary}
For an arbitrary  $x_0 \in S^{-}$ the following condition holds:
$$
x_0=\bigcap^{\infty} _{n=1} {\Delta^{-} _{c_1c_2...c_n}}.
$$
\end{corollary}

Let $u$  be a fixed positive integer from  $A$.

By $\mathbb S_{(-s,u)}$ denote the set (a subset of the segment $\left[- \frac{s}{s+1}, \frac{1}{s+1}\right]$) of all numbers $x$ represented by the nega-s-adic expansion such that are of the form
$$
 x=\sum^{\infty} _{n=1} {\left(\frac{\alpha_n -u }{{(-s)}^{\alpha_1+\alpha_2+\dots+\alpha_n}}\right)-\frac{u}{s+1}}, 
$$
where $(\alpha_n) \in L$.

This set is the following set
$$
\mathbb S_{(-s,u)}=\left\{x: x=\Delta^{-s} _{\underbrace{u\ldots u}_{{\alpha_1-1 }}\alpha_1 \underbrace{u\ldots u}_{\alpha_2-1}\alpha_2...\underbrace{u\ldots u}_{\alpha_n-1}\alpha_n...}\right\},
$$
where  $(\alpha_n)\in L$, $u \ne \alpha_n$ for all $n \in \mathbb N$, and $u$ is a fixed number. 

It was shown in \cite{S. Serbenyuk 2013} that the following statement is true.
\begin{theorem}
\label{the: osn(-s)}
Let  $\left\{\sigma_{1},\sigma_{2},\dots ,\sigma_{m}\right\}$ be a fixed finite set of combinations (tuples) of nega-s-adic digits, 
 $E$ be a set whose elements have in own nega-s-adic representation only combinations of digits from the set  $\left\{\sigma_{1},\sigma_{2},\dots ,\sigma_{m}\right\}$. Then  the Hausdorff-Besicovitch dimension $\alpha_0 (E)$  of the set   $E$  satisfies the equation
$$
N(\sigma^1 _m)\left(\frac{1}{s}\right)^{\alpha_0}+N(\sigma^2 _m)\left(\frac{1}{s}\right)^{2\alpha_0}+\dots+N(\sigma^{k} _m)\left(\frac{1}{s}\right)^{k\alpha_0}=1,
$$
where $N(\sigma^{k} _m)$ is a number of $k$-digit combinations from   $\left\{\sigma_{1},\sigma_{2},\dots ,\sigma_{m}\right\}$, 
$k \in \mathbb N$,  and $N(\sigma^{1} _m)+N(\sigma^{2} _m)+\dots+N(\sigma^{k} _m)=m$.
\end{theorem}
\begin{proof} Let $\left\{\sigma_{1},\sigma_{2},...,\sigma_{m}\right\}$ be a set of fixed combinations of nega-s-adic digits, and the nega-s-adic representation of any number from $E$ ($E$ is a  Cantor-like set) contains only such combinations of
digits. There exist digit combinations   $e_1e_2...e_r$, $\iota_1\iota_2...\iota_t$, where $r,t \in \mathbb N$ (they can be represented as one or several combinations from  $\left\{\sigma_{1},\sigma_{2},...,\sigma_{m}\right\}$) such that 
$$
\inf E= \Delta^{-s} _{(e_1e_2...e_r)(e_1e_2...e_r)...}~ \mbox{and}~ \sup E= \Delta^{-s} _{(\iota_1\iota_2...\iota_t)(\iota_1\iota_2...\iota_t)...}.
$$
Also, here
$$
d(E)=\sup E - \inf E, ~\mbox{where}~ d(\cdot) ~\mbox{is the diameter of the set}.
$$

A cylinder $\Delta^{(-s,E)} _{\tau_{1}\tau_{2}...\tau_{n}}$ of rank $n$ with the base  $\tau_{1}\tau_{2}...\tau_{n}$ is a set formed by all numbers
of $E$ with nega-s-adic representations in which the first $n$ combinations of digits are
fixed and are from $\left\{\sigma_{1},\sigma_{2},...,\sigma_{m}\right\}$. It is easy to see that 
$$
d\left(\Delta^{(-s,E)} _{\tau_{1}\tau_{2}...\tau_{n}}\right)=\frac{d(E)}{s^{N(\tau_{1}+\tau_{2}+...+\tau_{n})}},
$$
where $N(\tau_{1}+\tau_{2}+...+\tau_{n})$ is the number of digits in the combination  $\tau_{1}\tau_{2} ... \tau_{n}$.

Since
$$
E=C[-s,\left\{\sigma_{1},\sigma_{2},...,\sigma_{m}\right\} ], ~~~E\subset \left[\inf E; \sup E \right],~ \mbox{and}
$$
$$
\frac{\Delta^{(-s,E)} _{\tau_{1}\tau_{2}...\tau_{n}\tau_{n+1}}}{\Delta^{(-s,E)} _{\tau_{1}\tau_{2}...\tau_{n}}}=\frac{1}{s^{N(\tau_{n+1})}},
$$
we have
$$
E=\left[I_{\tau_1} \cap E\right]\cup\left[I_{\tau_2} \cap E\right]\cap...\cap \left[I_{\tau_m} \cap E\right],
$$
where $I_{\tau_i} =\left[\inf \Delta^{(-s,E)} _{\tau_{i}}; \sup \Delta^{(-s,E)} _{\tau_{i}} \right]$, $i=\overline{1,m}$.

So,
$$
\left[I_{\tau^1 _1}\cap E\right]\stackrel{s^{-1}}{\sim}E,  \left[I_{\tau^1 _2}\cap E\right]\stackrel{s^{-1}}{\sim}E,...,\left[I_{\tau^1 _{n_1}}\cap E\right]\stackrel{s^{-1}}{\sim}E;
$$
$$
\left[I_{\tau^2 _1}\cap E\right]\stackrel{s^{-2}}{\sim}E,  \left[I_{\tau^2 _2}\cap E\right]\stackrel{s^{-2}}{\sim}E,...,\left[I_{\tau^2 _{n_2}}\cap E\right]\stackrel{s^{-2}}{\sim}E;
$$
$$
................................................................
$$
$$
\left[I_{\tau^k _1}\cap E\right]\stackrel{s^{-k}}{\sim}E,  \left[I_{\tau^k _2}\cap E\right]\stackrel{s^{-k}}{\sim}E,...,\left[I_{\tau^k _{n_k}}\cap E\right]\stackrel{s^{-k}}{\sim}E,
$$
where $\tau^k _j$ is some    $k$-digit combination from  $\left\{\sigma_{1},\sigma_{2},...,\sigma_{m}\right\}$ ($j=\overline{1,{n_k}}$), and $n_k$ is the number of $k$-digit combinations from $\left\{\sigma_{1},\sigma_{2},...,\sigma_{m}\right\}$.

Hence  the set  $E$ is a self-similar fractal whose Hausdorff  dimension satisfies the equation 
$$
N(\sigma^1 _m)\left(\frac{1}{s}\right)^{\alpha_0}+N(\sigma^2 _m)\left(\frac{1}{s}\right)^{2\alpha_0}+...+N(\sigma^{k} _m)\left(\frac{1}{s}\right)^{k\alpha_0}=1.
$$
\end{proof}

The following statements follow from the last-mentioned theorem.
\begin{theorem}[\cite{S. Serbenyuk 2013}]
The set  $\mathbb S_{(-s,u)}$ is:
\begin{itemize}
\item an uncountable,   perfect,   nowhere dense sets of zero Lebesgue measure;
\item a self-similar fractal fractal, and its Hausdorff-Besicovitch dimension  $\alpha_0(\mathbb S_{(-s,u)})$ satisfies the equation
$$
\sum_{i \in A_u} {\left(\frac{1}{s}\right)^{i\alpha_0}}=1, ~\mbox{where} ~A_u= \{1,2,\dots ,s-1\}\setminus\{ u\}.
$$
\end{itemize}
\end{theorem}

Let us consider fractal sets whose elemets represented by nega-s-adic Cantor series. 
\begin{theorem}[\cite{S. Serbenyuk 2013}] 
Let $s>1$ be a fixed positive integer, $\alpha_{m_1+m_2+\dots+m_n}\ne 0$ for all $n \in \mathbb N$, and   $m_n\in \{3,5,7,\dots,2i+1,\dots\}$. Then the set  $M_{(-D,s)}$  
$$
M_{(-D,s)}=\left\{x: x= \Delta^{-s} _{\underbrace{0\ldots 0}_{m_1-1}\alpha_{m_1}\underbrace{0\ldots 0}_{m_2-1}\alpha_{m_1+m_2}\ldots\underbrace{0\ldots 0}_{m_n-1}\alpha_{m_1+m_2+\dots+m_n}...}\right\}.
$$
is a self-similar fractal whose Hausdorff-Besicovitch dimension $\alpha_0\left(M_{(-D,s)}\right)$ is equal to 
$$
\log_s {\left(\sqrt[3]{\frac{s-1}{2}+\frac{1}{6}\sqrt{\frac{27(s-1)^2-4}{3}}}+\sqrt[3]{\frac{s-1}{2}-\frac{1}{6}\sqrt{\frac{27(s-1)^2-4}{3}}}\right)}.
$$
\end{theorem}
\begin{proof}
From \eqref{nsc1} and Theorem  \ref{the: osn(-s)}, it follows that the Hausdorff-Besicovitch dimension of the set  $M_{(-D,s)}$ under $m_n \in\{3,5,7,...,2i+1,...\}$, $\alpha_{m_1+m_2+...+m_n}\ne 0$, and under  fixed $s>1$, satisfies the equation
$$
(s-1)\left(\frac{1}{s}\right)^{3\alpha_0}+(s-1)\left(\frac{1}{s}\right)^{5\alpha_0}+(s-1)\left(\frac{1}{s}\right)^{7\alpha_0}+...+(s-1)\left(\frac{1}{s}\right)^{(2i+1)\alpha_0}+...=1,~i=1,2,....
$$
The last equation is equivalent to the equation  
$$
s^{3\alpha_0}-s^{\alpha_0}-(s-1)=0.
$$
Using Cardano's formula, we get the result.
\end{proof}
%\begin{remark}
%ßêùî óìîâà $m_n=1 $  âèêîíóºòüñÿ ñê³í÷åííó ê³ëüê³ñòü ðàç³â, òî çíà÷åííÿ ðîçì³ðíîñò³ ìíîæèíè $M_{(-D,s)}$ íå %çì³íþºòüñÿ. ßêùî æ óìîâà  $m_n=1$ âèêîíóºòüñÿ äëÿ âñ³õ~$n \in~\mathbb N$   àáî $m_n \ne 1$  âèêîíóºòüñÿ ñê³í÷åííó %ê³ëüê³ñòü ðàç³â, òîä³ çíà÷åííÿ ðîçì³ðíîñò³ $\alpha_0(M_{(-D,s)})$ Õàóñäîðôà-Áåçèêîâè÷à ìíîæèíè $M_{(-D,s)}$ %äîð³âíþâàòèìå  $\log_s {(s-1)}$.
%\end{remark}

\begin{corollary}
If a sequence $(m_n)$  of odd positive integers is a fixed  purely periodic sequence with the period  $(m_1m_2\ldots m_t)$, then 
the set $M^{'} _{(-D,s,t)}$ of all numbers represented by nega-s-adic Cantor series \eqref{nsc1} with the corresponding sequence 
$(m_n)$ is a self-similar fractal and 
$$
\alpha_0\left(M^{'} _{(-D,s,t)}\right)=\frac{t}{m_1+m_2+\dots+m_t}.
$$
\end{corollary}
\begin{proof} Since elements of this set have periodic nega-s-adic representation, i.e., 
$$
M^{'} _{(-D,s,t)} \ni x=\Delta^{-s} _{\left(\underbrace{0...0}_{m_1-1}\alpha_{m_1}\underbrace{0...0}_{m_2-1}\alpha_{m_1+m_2}...\underbrace{0...0}_{m_t-1}\alpha_{m_1+m_2+...+m_t}\right)},
$$
where $\{m_1,m_2,...,m_t\}$ is a fixed set of odd numbers and $\alpha_{m_1},\alpha_{m_1+m_2},...,\alpha_{m_1+...+m_t}$ are numbers from  the set $A$, from Theorem 
 \ref{the: osn(-s)}, it follows that Hausdorff-Besicovitch  dimension satisfies the equation 
$$
s^t\left(\frac{1}{s}\right)^{(m_1+m_2+...+m_t)\alpha_0}=1.
$$
The statement follows from the last equation.
\end{proof}

 Finally, let us remark that restrictions on using elements of sets $\mathbb S_{(\pm s,u)}$ are new (they occur for the first time).

So, we considered topological, metric, and fractal properties of certain sets whose elements have restrictions on using digits in own expansions. For considered sets, the case of functional restrictions is equivalent to the case of restrictions on using combinations of digits. The simple methods for the calculation of the Hausdorff-Besicovitch dimension of such sets are described. 
In the  case of the s-adic or nega-s-adic representations, the Hausdorff-Besicovitch dimension of a set whose elements have in own representations only combinations of digits from some fixed set of combinations of digits, depends on parameters as a number of $k$-digit combinations and numbers $k$. 
 In addition, note that considered sets have the Moran structure. Similar investigations did not study for the case of generalizations of the s-adic or nega-s-adic representation. These investigations will be discussed  by the author of the present article in a further paper.

%%%%%%%%%%%%%%%%%%%%%%%%
\section{Properties of images}
%%%%%%%%%%%%%%%%%%%%%%%%

In this section, the main attention is given to images of sets  $\mathbb S_{(s,u)}$ and  $\mathbb S_{(- s,u)}$ under the Salem type functions   (see \cite{Symon2015, Symon2017, Symon2019}, the Salem function was introduced in \cite{Salem1943}). 

Let $s>1$ be a fixed positive integer and $\alpha_n\in A=\{0,1,\dots, s-1\}$. Let  $P=\{p_0,p_1,\dots , p_{s-1}\}$ be a fixed set whose elements satisfy the following properties: $p_0+p_1+\dots+p_{s-1}=1$ and $p_i>0$ for all $i=\overline{0,s-1}$. Then let us consider the following distribution functions. 

Let $\zeta$ be a random variable defined by the s-adic representation
$$
\zeta= \frac{\iota_1}{s}+\frac{\iota_2}{s^2}+\frac{\iota_3}{s^3}+\dots +\frac{\iota_{k}}{s^{k}}+\dots   = \Delta^{s} _{\iota_1\iota_2...\iota_{k}...},
$$
where  digits $\iota_k$ $(k=1,2,3, \dots )$ are random and taking the values $0,1,\dots ,s-1$ with positive probabilities ${p}_{0}, {p}_{1}, \dots , {p}_{s-1}$. That is,  $\iota_k$ are independent and  $P\{\iota_k=\alpha_k\}={p}_{\alpha_k}$, $\alpha_k \in A$. 

Let $\varsigma$ be a random variable defined by the s-adic representation
$$
\varsigma= \Delta^{s} _{\pi_1\pi_2...\pi_{k}...}=\sum^{\infty} _{k=1}{\frac{\pi_k}{s^k}},
$$
where 
$$
\pi_k=\begin{cases}
\alpha_k&\text{if  $k$ is odd}\\
s-1-\alpha_k&\text{if $k$ is even}
\end{cases}
$$
and digits $\pi_k$ $(k=1,2,3, \dots )$ are random and taking the values $0,1,\dots ,s-1$ with positive probabilities ${p}_{0}, {p}_{1}, \dots , {p}_{s-1}$. That is, $\pi_k$ are independent and  $P\{\pi_k=\alpha_k\}={p}_{\alpha_k}$, $P\{\pi_k=s-1-\alpha_k\}={p}_{s-1-\alpha_k}$, where $\alpha_k \in A$.

Let us consider the distribution function  ${f}_{\zeta}$ of the random variable $\zeta$ and the distribution function  ${\tilde F}_{\varsigma}$ of the random variable $\varsigma$:
$$
{f}_{\zeta}(x)=\begin{cases}
0&\text{whenever $x< 0$}\\
\beta_{\alpha_1(x)}+\sum^{\infty} _{k=2} {\left({\beta}_{\alpha_k(x)} \prod^{k-1} _{j=1} {{p}_{\alpha_j(x)}}\right)}&\text{whenever $0 \le x<1$}\\
1&\text{whenever $x\ge 1$,}
\end{cases}
$$
where ${p}_{\alpha_{j(x)}}>0$ and
$$
\beta_{\alpha_k}=\begin{cases}
\sum^{\alpha_k(x)-1} _{i=0} {p_{i}(x)}&\text{whenever $\alpha_k(x)>0$}\\
0&\text{whenever  $\alpha_k(x)=0$;}
\end{cases}
$$
also,
$$
{\tilde F}_{\varsigma}(x)=\begin{cases}
0&\text{whenever $x< 0$}\\
\tilde\beta_{\alpha_1(x)}+\sum^{\infty} _{k=2} {\left({\tilde\beta}_{\alpha_k(x)} \prod^{k-1} _{j=1} {{\tilde p}_{\alpha_j(x)}}\right)}&\text{whenever $0 \le x<1$}\\
1&\text{whenever $x\ge 1$,}
\end{cases}
$$
where ${p}_{\alpha_{j(x)}}>0$,
$$
x=\Delta^s _{\tilde\alpha_1\tilde\alpha_2...\tilde\alpha_k...}=\frac{1}{s+1}-\Delta^{-s} _{\alpha_1\alpha_2...\alpha_k...}\equiv\frac{1}{s+1}-\sum^{\infty} _{k=1}{\frac{(-1)^k\alpha_{k}}{s^{k}}}=\sum^{\infty} _{k=1}{\frac{\alpha_{2k-1}}{s^{2k-1}}}+\sum^{\infty} _{k=1}{\frac{s-1-\alpha_{2k}}{s^{2k}}},
$$
and 
$$
\tilde  p_{\alpha_k}=\begin{cases}
p_{\alpha_k}&\text{if  $k$ is odd}\\
p_{s-1-\alpha_k}&\text{if $k$ is even},
\end{cases}
$$
$$
\tilde\beta_{\alpha_k}=\begin{cases}
\beta_{\alpha_k}&\text{if  $k$ is odd}\\
\beta_{s-1-\alpha_k}&\text{if $k$ is even}.
\end{cases}
$$

One can note that  the function  
$$ 
{\tilde F}(x)=\beta_{\alpha_1(x)}+\sum^{\infty} _{n=2} {\left({\tilde\beta}_{\alpha_n(x)}\prod^{n-1} _{j=1} {{\tilde p}_{\alpha_j(x)}}\right)},
$$
is a partial case of the function investigated in \cite{Symon2019}. 

Let $x\in\mathbb S_{(s,u)}$. Let us consider properties of the following  images of  $\mathbb S_{(-s,u)}$ and $\mathbb S_{(s,u)}$:
$$
\mathbb S_{(P,u)}=\{y: y=f_{\xi}(x), x\in \mathbb S_{(s,u)}\}
$$
and
$$
\mathbb S_{(-P,u)}=\left\{\tilde y: \tilde y=\tilde F \circ f_l\circ f_+(x), x\in\mathbb S_{(s,u)} \right\}=\{z: z=\tilde F\circ f_l (x),~~~ x\in \mathbb S_{(-s,u)}\}.
$$ 
Here
$$
\tilde y=\tilde F \circ f_l\circ f_+(x),
$$
where 
$$
f_+: x=\Delta^s _{\alpha_1\alpha_2...\alpha_n...}\to \Delta^{-s} _{\alpha_1\alpha_2...\alpha_n...}=y
$$
 is not monotonic on the domain and is a nowhere differentiable function~(\cite{S. Serbenyuk functions with complicated local structure 2013}), $f_l(y)=\frac{1}{s+1}-y$, and $\tilde F$ is the last-mentioned distribution function.

Let us describe  properties of the set $\mathbb S_{(P,u)}$.
\begin{theorem}{\cite{Symon}.}
\label{th: S(p,u)}
The set $\mathbb S_{(P,u)}$ is an uncountable, perfect, and nowhere dense set of zero Lebesgue measure and also is  a self-similar fractal whose Hausdorff dimension $\alpha_0 (\mathbb S_{(P,u)})$  satisfies the following equation 
$$
\sum _{i\in A_0\setminus\{u\}} {\left(p_ip^{i-1} _u\right)^{\alpha_0}}=1.
$$
\end{theorem}

Let $c_1, c_2,\dots , c_n$ be an ordered tuple of integers such that $c_i\in\{0,1,\dots ,s-1\}$ for $i=\overline{1,n}$.
\begin{definition} 
{\itshape A cylinder of rank $n$ with  base $c_1c_2\ldots c_n$} is a set $\Delta^{(P,u)} _{c_1c_2\ldots c_n}$ of the form: 
$$
\Delta^{(P,u)} _{c_1c_2\ldots c_n}=\left\{x: x=\Delta^{P}_{{\underbrace{u...u}_{c_1-1}} c_1{\underbrace{u...u}_{c_2 -1}}c_2 ...{\underbrace{u...u}_{ c_n -1}}c_n{\underbrace{u...u}_{\alpha_{n+1}-1}}\alpha_{n+1}{\underbrace{u...u}_{\alpha_{n+2}-1}}\alpha_{n+2}...}, \alpha_j=c_j, j=\overline{1,n}\right\}.
$$
\end{definition}

By $(a_1a_2\ldots a_k)$ denote the period $a_1a_2\ldots a_k$ in the representation of a periodic number.

The following lemma describes local properties of the set $\mathbb S_{(P,u)}$.
\begin{lemma}{\cite{Symon}.} Cylinders $ \Delta^{(P,u)} _{c_1...c_n} $ have the following properties:
\label{lm: Lemma on cylinders}
\begin{enumerate}
\item
$$
\inf  \Delta^{(P,u)} _{c_1...c_n}=\begin{cases}
\Delta^{P} _{{\underbrace{0...0}_{c_1-1}} c_1{\underbrace{0...0}_{c_2 -1}}c_2 ...{\underbrace{0...0}_{ c_n -1}}c_n({\underbrace{0...0}_{ s-2}}[s-1])} &\text{if $u=0$}\\
\Delta^{P} _{{\underbrace{1...1}_{c_1-1}} c_1{\underbrace{1...1}_{c_2 -1}}c_2 ...{\underbrace{1...1}_{ c_n -1}}c_n({\underbrace{1...1}_{ s-2}}[s-1])} &\text{if $u=1$}\\
$$\\
\Delta^{P} _{{\underbrace{u...u}_{c_1-1}} c_1{\underbrace{u...u}_{c_2 -1}}c_2 ...{\underbrace{u...u}_{ c_n -1}}c_n(1)}&\text{if $ u \in \{2,3,\dots ,s-1\}$,}
\end{cases}
$$
$$
\sup  \Delta^{(P,u)} _{c_1...c_n...}=\begin{cases}
\Delta^{P} _{{\underbrace{[s-1]...[s-1]}_{c_1-1}} c_1 ...{\underbrace{[s-1]...[s-1]}_{ c_n -1}}c_n({\underbrace{[s-1]...[s-1]}_{ s-3}}[s-2])} &\text{if $u=s-1$}\\
\Delta^{P} _{{\underbrace{u...u}_{c_1-1}} c_1{\underbrace{u...u}_{c_2 -1}}c_2 ...{\underbrace{u...u}_{ c_n -1}}c_n({\underbrace{u...u}_{ u}}[u+1])} &\text{if $u\in\{1,\dots, s-2\}$}\\
$$\\
\Delta^{P} _{{\underbrace{0...0}_{c_1-1}} c_1{\underbrace{0...0}_{c_2 -1}}c_2 ...{\underbrace{0...0}_{ c_n -1}}c_n(1)}&\text{if $ u=0$.}
\end{cases}
$$
\item If $d(\cdot) $ is the diameter of a set, then
$$
d(\Delta^{(P,u)} _{c_1...c_n})=d(\mathbb S_{(P_s,u)})p^{c_1+c_2+\dots+c_n-n} _{u}\prod^{n} _{j=1}{p_{c_j}}.
$$
\item
$$
\frac{d(\Delta^{(P,u)} _{c_1...c_nc_{n+1}})}{d(\Delta^{(P,u)} _{c_1...c_n})}=p_{c_{n+1}}p^{c_{n+1}-1} _{u}.
$$
\item 
$$
  \Delta^{(P,u)} _{c_1c_2...c_n} =\bigcup^{s-1} _{i=1} { \Delta^{(P,u)} _{c_1c_2...c_ni}}~~~\forall c_n \in A_0,~~~n \in \mathbb N,~ i \ne u.
$$
\item The following relationships are satisfied: 
\begin{enumerate}
\item if $ u\in \{0,1\}$, then 
$$
\inf \Delta^{(P,u)} _{c_1...c_np}> \sup \Delta^{(P,u)} _{c_1...c_n[p+1]};
$$
\item if  $ u \in \{2,3,\dots ,s-3\}$, then 
$$
\begin{cases}
\sup \Delta^{(P,u)} _{c_1...c_np}< \inf \Delta^{(P,u)} _{c_1...c_n[p+1]}&\text{for all $p+1\le u$}\\
$$\\
\inf \Delta^{(P,u)} _{c_1...c_np}> \sup \Delta^{(P,u)} _{c_1...c_n[p+1]},&\text{for all $u<p$;}
\end{cases}
$$
\item if $ u  \in \{s-2,s-1\}$, then
$$
\sup \Delta^{(P,u)} _{c_1...c_np}< \inf \Delta^{(P,u)} _{c_1...c_n[p+1]} ~~~(\text{in this case, the condition $p\ne s-1$ holds}).
$$
\end{enumerate}
\end{enumerate}
\end{lemma}

We considered properties of $\mathbb S_{(s,u)}$ and its image $\mathbb S_{(P,u)}$  under the Salem function. So, the Salem function preserves the self-similarity, but, in the general case, does not preserve the Hausdorff dimension. This map also  preserves the structure of   $\mathbb S_{(s,u)}$ but numerical values change.

Finally, one can note the following theorem.
\begin{theorem}{\cite{Symon}.}
Let $\mathbb S$ be a set whose elements represented in terms of the s-adic representation by a finite number of fixed combinations $\tau_1, \tau_2,\dots,\tau_m$ of  digits from the alphabet $A$.

Let  $E$ be an image of the set $\mathbb S$ under the Salem function $f_{\xi}$.  Then the Hausdorff  dimension $\alpha_0$ of $E$ satisfies the following equation: 
$$
\sum^{m} _{j=1}{\left(\prod^{s-1} _{i=0}{p^{N_i(\tau_j)} _i}\right)^{\alpha_0}}=1,
$$
where $N_i(\tau_k)$ ($k=\overline{1,m}$)  is a number of  the digit $i$ in $\tau_k$ from the set $\{\tau_1, \tau_2,\dots,\tau_m\}$.
\end{theorem}

Now we describe properties of $\mathbb S_{(-P,u)}$. 

Suppose $d(\cdot)$ is the diameter of a set and a cylinder $\Delta^{(-P,u)} _{c_1c_2...c_n}$ is a set whose elements are elements of $\mathbb S_{(-P,u)}$ and for these elements the condition $\alpha_i=c_i$ holds for all $i=\overline{1,n}$ (here $c_1, c_2,\cdots, c_n$ is a fixed tuple).  

\begin{theorem}{\cite{Serbenyuk2021}.}
\label{th: the first main theorem}
An arbitrary set $\mathbb S_{(-P,u)}$ is an uncountable, perfect, and nowhere dense set of zero Lebesgue measure.
\end{theorem}

\begin{theorem}{\cite{Serbenyuk2021}.}
\label{th: the second main theorem}
In the general case,  the set  $\mathbb S_{(-P,u)} $ is not a self-similar fractal,  the Hausdorff dimension $\alpha_0 (\mathbb S_{(-P,u)})$ of which can be calculated by the formula: 
$$
\alpha_0=\liminf_{k\to\infty}{\alpha_k},
$$
where $(\alpha_k)$ is a sequence of numbers satisfying the equation
$$
\left(\sum_{\substack{c_1 \text{is odd}\\ c_1\in \overline{A}}}{\left(\omega_{2,c_1}\right)^{\alpha_1} }+\sum_{\substack{c_1 \text{is even}\\ c_1\in \overline{A}}}{\left(\omega_{4,c_1}\right)^{\alpha_1} }\right)\times
$$
$$
\times\prod^k _{i=2}{\left(\sum_{\substack{c_i \text{is odd}\\ c_i\in \overline{A}}}{N_{1,c_i}\left(\omega_{1,c_i}\right)^{\alpha_i} }+\sum_{\substack{c_i \text{is odd}\\ c_i\in \overline{A}}}{N_{2,c_i}\left(\omega_{2,c_i}\right)^{\alpha_i} }+\sum_{\substack{c_i \text{is even}\\ c_i\in \overline{A}}}{N_{3,c_i}\left(\omega_{3,c_i}\right)^{\alpha_i} }+\sum_{\substack{c_i \text{is even}\\ c_i\in \overline{A}}}{N_{4,c_i}\left(\omega_{4,c_i}\right)^{\alpha_i} }\right)} =1.   
$$
Here $N_{j,c_i}$ ($j=\overline{1,4}, 1<i\in\mathbb N$) is the number of cylinders $\Delta^{(-P,u)} _{c_1c_2...c_i}$ for which
$$
\frac{d\left(\Delta^{(-P,u)} _{c_1c_2...c_{i-1}c_i}\right)}{d\left(\Delta^{(-P,u)} _{c_1c_2...c_{i-1}}\right)}=\omega_{j,c_i}.
$$
Also, 
$$
\omega_{1,c_i}=\underbrace{p_{s-1-u}p_u\ldots p_{s-1-u}p_u}_{c_i-1}p_{s-1-c_i}\frac{d\left(\overline{\mathbb S_{(P,u)}}\right)}{d\left(\underline{\mathbb S_{(P,u)}}\right)}~~~\text{for an  odd number $c_i$},
$$
$$
\omega_{2,c_i}=\underbrace{p_up_{s-1-u}\ldots p_up_{s-1-u}}_{c_i-1}p_{c_i}\frac{d\left(\underline{\mathbb S_{(P,u)}}\right)}{d\left(\overline{\mathbb S_{(P,u)}}\right)}~~~\text{for an odd number $c_i$},
$$
$$
\omega_{3,c_i}=\underbrace{p_{s-1-u}p_u\ldots p_{s-1-u}p_up_{s-1-u}}_{c_i-1}p_{c_i}~~~\text{for an even number $c_i$},
$$
$$
\omega_{4,c_i}=\underbrace{p_up_{s-1-u}\ldots p_up_{s-1-u}p_{u}}_{c_i-1}p_{s-1-c_i}~~~\text{for an even number $c_i$}.
$$
In addition, $N_{1,c_i}+N_{2,c_i}=l(m+l)^{i-1}$ and $N_{3,c_i}+N_{4,c_i}=m(m+l)^{i-1}$, where $l$ is the number of odd numbers in the set  $\overline{A}=A\setminus\{0,u\}$ and $m$ is the number of even numbers in $\overline{A}$.
\end{theorem}
Auxiliary values can be calculated from the following lemma.
\begin{lemma}{\cite{Serbenyuk2021}.}
\label{lm: auxiliary lemma 6}
For the sets $\underline{\mathbb S_{(P,u)}}$ and $\overline{\mathbb S_{(P,u)}}$, the following equalities hold: 
$$
\inf {\underline{\mathbb S_{(P,u)}}}=\begin{cases}
\Delta^{P} _{[s-2](0[s-3])}&\text{if $u=0$}\\
\Delta^{P} _{[s-2]1[s-4](1[s-3])}&\text{if $u=1$}\\
\Delta^{P} _{([s-1-u]2)}&\text{if $u\in\{2,3,\dots , s-1\}$}
\end{cases}
$$
and
$$
\sup {\underline{\mathbb S_{(P,u)}}}=\begin{cases}
\Delta^{P} _{([s-1-u]2)}&\text{if $u\in\{0,1\}$}\\
\Delta^{P} _{[s-2](u[s-3])}&\text{if $u\in\{2,3,\dots , s-1\}$},
\end{cases}
$$
$$
 \inf {\overline{\mathbb S_{(P,u)}}}=\begin{cases}
\Delta^{P} _{(u[s-3])}&\text{if $u\in\{0,1\}$}\\
\Delta^{P} _{1([s-1-u]2)}&\text{if $u\in\{2,3,\dots , s-1\}$}
\end{cases}
$$
and
$$
 \sup {\overline{\mathbb S_{(P,u)}}}=\begin{cases}
\Delta^{P} _{1([s-1]2)}&\text{if $u=0$}\\
\Delta^{P} _{1[s-2]3([s-2]2)}&\text{if $u=1$}\\
\Delta^{P} _{(u[s-3])}&\text{if $u\in\{2,3,\dots , s-1\}$}
\end{cases}
$$
\end{lemma}

Finally, let us consider local properties of $\mathbb S_{(-P,u)} $. 

Assume  
$$
\tilde \alpha_n=\begin{cases}
\alpha_n &\text{whenever $n$ is odd}\\
s-1-\alpha_n &\text{whenever $n$ is even}
\end{cases}
$$
and
$$
\tilde u=\begin{cases}
u &\text{whenever $u$ is situated  at an odd position in the representation}\\
s-1-u &\text{whenever $u$ is situated at  in an even position in the representation}
\end{cases}.
$$

\begin{lemma} Cylinders $ \Delta^{(-P,u)} _{c_1...c_n} $ have the following properties:
\label{lm: Lemma on cylinders}
\begin{enumerate}
\item
$$
\inf  \Delta^{(-P,u)} _{c_1...c_n}=\begin{cases}
\tau_n+\left(\prod^{n} _{j=1}{\tilde p_{c_j, c_1+...+c_j}}\right)\left(\prod_{\substack{i=\overline{1,c_1+...+c_n-1}\\ i\notin C_{n-1}}}{\tilde p_{u,i}}\right)\inf{\overline{\mathbb S_{(P,u)}}} &\text{if $c_1+\dots +c_n$ is even} \\
\tau_n+\left(\prod^{n} _{j=1}{\tilde p_{c_j, c_1+...+c_j}}\right)\left(\prod_{\substack{i=\overline{1,c_1+...+c_n-1}\\ i\notin C_{n-1}}}{\tilde p_{u,i}}\right)\inf{\underline{\mathbb S_{(P,u)}}} &\text{if $c_1+\dots +c_n$ is odd} \\
\end{cases},
$$
$$
\sup  \Delta^{(-P,u)} _{c_1...c_n}=\begin{cases}
\tau_n+\left(\prod^{n} _{j=1}{\tilde p_{c_j, c_1+...+c_j}}\right)\left(\prod_{\substack{i=\overline{1,c_1+...+c_n-1}\\ i\notin C_{n-1}}}{\tilde p_{u,i}}\right)\sup{\overline{\mathbb S_{(P,u)}}} &\text{if $c_1+\dots +c_n$ is even} \\
\tau_n+\left(\prod^{n} _{j=1}{\tilde p_{c_j, c_1+...+c_j}}\right)\left(\prod_{\substack{i=\overline{1,c_1+...+c_n-1}\\ i\notin C_{n-1}}}{\tilde p_{u,i}}\right)\sup{\underline{\mathbb S_{(P,u)}}} &\text{if $c_1+\dots +c_n$ is odd} \\
\end{cases},
$$
where 
$$
\tau_n=\Delta^P _{{\underbrace{\tilde u...\tilde u}_{c_1-1}} \tilde c_1{\underbrace{\tilde u...\tilde u}_{c_2 -1}}\tilde c_2 ...{\underbrace{\tilde u...\tilde u}_{ c_n -1}}\tilde c_n(0)}.
$$

\item If $d(\cdot) $ is the diameter of a set, then
$$
d\left( \Delta^{(-P,u)} _{c_1...c_n}\right)=\begin{cases}
\left(\prod^{n} _{j=1}{\tilde p_{c_j, c_1+...+c_j}}\right)\left(\prod_{\substack{i=\overline{1,c_1+...+c_n-1}\\ i\notin C_{n-1}}}{\tilde p_{u,i}}\right)d\left({\overline{\mathbb S_{(P,u)}}}\right) &\text{if $c_1+\dots +c_n$ is even} \\
\left(\prod^{n} _{j=1}{\tilde p_{c_j, c_1+...+c_j}}\right)\left(\prod_{\substack{i=\overline{1,c_1+...+c_n-1}\\ i\notin C_{n-1}}}{\tilde p_{u,i}}\right)d\left({\underline{\mathbb S_{(P,u)}}}\right) &\text{if $c_1+\dots +c_n$ is odd} \\
\end{cases}
$$
\item $\frac{d(\Delta^{(-P,u)} _{c_1...c_nc_{n+1}})}{d(\Delta^{(-P,u)} _{c_1...c_n})}=$
$$
=\begin{cases}
p_{s-1-c_{n+1}}\left(\prod^{c_1+...+c_{n+1}-1} _{{i={c_1+c_2+...+c_n+1}}}{\tilde p_{u,i}}\right)&\text{if $c_1+\dots +c_n$, $c_{n+1}$ are  even} \\
p_{c_{n+1}}\left(\prod^{c_1+...+c_{n+1}-1} _{{i={c_1+c_2+...+c_n+1}}}{\tilde p_{u,i}}\right)&\text{if $c_1+\dots +c_n$ is odd,  $c_{n+1}$ is  even} \\
p_{s-1-c_{n+1}}\left(\prod^{c_1+...+c_{n+1}-1} _{{i={c_1+c_2+...+c_n+1}}}{\tilde p_{u,i}}\right)\left({\overline{\mathbb S_{(P,u)}}}/{\underline{\mathbb S_{(P,u)}}}\right) &\text{if $c_1+\dots +c_n$, $c_{n+1}$ are odd } \\
p_{c_{n+1}}\left(\prod^{c_1+...+c_{n+1}-1} _{{i={c_1+c_2+...+c_n+1}}}{\tilde p_{u,i}}\right)\left({\underline{\mathbb S_{(P,u)}}}/{\overline{\mathbb S_{(P,u)}}}\right) &\text{if $c_1+\dots +c_n$ is even, $c_{n+1}$ is odd } \\
\end{cases}.
$$
\item 
$$
  \Delta^{(-P,u)} _{c_1c_2...c_n} =\bigcup_{c\in\overline{A}} { \Delta^{(-P,u)} _{c_1c_2...c_nc}}~~~\forall c_n \in \overline{A},~~~n \in \mathbb N.
$$
\item The following relationships are satisfied: 
\begin{enumerate}
\item if $ u\in \{0,1\}$, then 
$$
\begin{cases}
\inf \Delta^{(-P,u)} _{c_1...c_n[c+1]}> \sup \Delta^{(-P,u)} _{c_1...c_nc}&\text{whenever $c_1+\dots +c_n+c$ is even}\\
$$\\
\inf \Delta^{(-P,u)} _{c_1...c_nc}> \sup \Delta^{(-P,u)} _{c_1...c_n[c+1]}&\text{whenever $c_1+\dots +c_n+c$ is odd}
\end{cases}~~~(c\ne s-1);
$$
\item if  $ u \in \{2,3,\dots ,s-3\}$, then for an odd $c_1+\dots +c_n+c$ 
$$
\begin{cases}
\sup \Delta^{(-P,u)} _{c_1...c_nc}< \inf \Delta^{(P,u)} _{c_1...c_n[c+1]}&\text{for all $c+1\le u$}\\
$$\\
\inf \Delta^{(-P,u)} _{c_1...c_nc}> \sup \Delta^{(-P,u)} _{c_1...c_n[c+1]},&\text{for all $u<c$;}
\end{cases}
$$
if  $ u \in \{2,3,\dots ,s-3\}$, then for an even $c_1+\dots +c_n+c$ 
$$
\begin{cases}
\inf \Delta^{(-P,u)} _{c_1...c_n[c+1]}> \sup \Delta^{(-P,u)} _{c_1...c_nc}&\text{for all $u<c$}\\
$$\\
\inf \Delta^{(-P,u)} _{c_1...c_nc}> \sup \Delta^{(-P,u)} _{c_1...c_n[c+1]}&\text{for all $c+1\le u$}
\end{cases};
$$

\item if $ u  \in \{s-2,s-1\}$, then
$$
\begin{cases}
\inf \Delta^{(-P,u)} _{c_1...c_n[c+1]}> \sup \Delta^{(-P,u)} _{c_1...c_nc}&\text{whenever $c_1+\dots +c_n+c$ is odd}\\
\inf \Delta^{(-P,u)} _{c_1...c_nc}> \sup \Delta^{(-P,u)} _{c_1...c_n[c+1]}&\text{whenever $c_1+\dots +c_n+c$ is even}
\end{cases}.
$$
\end{enumerate}
\end{enumerate}
\end{lemma}

\begin{remark}
One can note that  if for $\mathbb S_{(-s,u)}$ and $\mathbb S_{(s,u)}$,  topological,  metric, and  fractal properties (without some properties of cylinders) are similar, then   fractal  and some local  properties of $\mathbb S_{(-P,u)}$ and $\mathbb S_{(P,u)}$ are different. For example, $\mathbb S_{(P,u)}$ is a self-similar fractal (i.e., this is a Moran set by Moran's definition, \cite{Moran1946}) but $\mathbb S_{(-P,u)}$ is  a non-self-similar set having the Moran structure (i.e., this is a  Moran set by the definition of Hua et al. (see the  definition in \cite{HRW2000})). 
\end{remark}

%%%%%%%%%%%%%%%%%%%%%%%%%
\section{Certain examples}

%%%%%%%%%%%%%%%%%%%%%%%%%%

Let us consider the case of the sets $\mathbb S_{(P_3,0)}$ and $\mathbb S_{(-P_3,0)}$. That is, the set $\mathbb S_{(P_3,0)}$ is a set of the form 
$$
 S_{(P_3,0)}:= \left\{x: x=\Delta^{P_3}_{{\underbrace{0...0}_{\alpha_1-1}}\alpha_1{\underbrace{0...0}_{\alpha_2 -1}}\alpha_2 ...{\underbrace{0...0}_{ \alpha_n -1}}\alpha_n...}, \alpha_n \in \{1,2\}, n=1,2,3, \dots \right\}.
$$
In the other words, our set is the image of a certain set under the Salem function $f_{\xi}$, and this certain set is the set whose elements represented in terms of the 3-adic (ternary) representation by using  combinations of ternary digits only from $\{1, 02\}$. So, applying Theorem~\ref{th: S(p,u)}, we obtain that the Hausdorff dimension of $\mathbb S_{(P_3,0)}$ satisfies the  equation
\begin{equation}
\label{eq: example-eq-1}
(p_1)^{\alpha_0}+(p_0p_2)^{\alpha_0}=1.
\end{equation}

Let us consider the set  $\mathbb S_{(-P_3,0)}$. That is, 
$$
 S_{(-P_3,0)}:= \left\{x: x=\Delta^{-P_3}_{{\underbrace{0...0}_{\alpha_1-1}}\alpha_1{\underbrace{0...0}_{\alpha_2 -1}}\alpha_2 ...{\underbrace{0...0}_{ \alpha_n -1}}\alpha_n...}, \alpha_n \in \{1,2\}, n=1,2,3, \dots \right\}.
$$
Since (\cite{Serbenyuk2021, 2016Serbeyuk-KNU}) 
$$
\Delta^{P_3} _{\alpha_1\alpha_2...\alpha_n...}\equiv \Delta^{-P_3} _{\alpha_1[2-\alpha_2]\alpha_3...\alpha_{2k-1}[2-\alpha_{2k}]...}
$$
and
$$
\Delta^{-P_3} _{\alpha_1\alpha_2...\alpha_n...}\equiv \Delta^{P_3} _{\alpha_1[2-\alpha_2]\alpha_3...\alpha_{2k-1}[2-\alpha_{2k}]...},
$$
we have that this set is a subset of the set
$$
\left\{x: x=\Delta^{P_3} _{\delta_1\delta_2...\delta_n...}, \delta_n \in\{00, 1, 22\}\right\}.
$$
Really, it is a subset, because, for example, 
$$
1=\Delta^{P_3} _{2222222...}\notin  \mathbb S_{(-P_3,0)}.
$$

Using Lemma~\ref{lm: auxiliary lemma 6}, we have
$$
\inf {{\mathbb S_{(-P_3,0)}}}=\inf {\overline{\mathbb S_{(P_3,0)}}}=\Delta^{-P_3} _{020202...}=\Delta^{P_3} _{0000000...}=0,
 $$
$$
\sup {{\mathbb S_{(-P_3,0)}}}=\sup {\overline{\mathbb S_{(P_3,0)}}}=\Delta^{-P_3} _{1020202...}=\Delta^{P_3} _{12222222...}=p_0+p_1 
$$
and
$$
\inf {\underline{\mathbb S_{(P_3,0)}}}=\Delta^{P_3} _{1000000...}=\beta_1=p_0,
$$
$$
\sup {\underline{\mathbb S_{(P_3,0)}}}=\Delta^{P_3} _{2222222...}=1. 
$$
Hence, using~Theorem~\ref{th: the second main theorem} and Lemma~\ref{lm: Lemma on cylinders}, we obtain (here $d(\cdot)=\sup(\cdot)-\inf(\cdot)$)
$$
d(\overline{\mathbb S_{(P_3,0)}})=\beta_2=p_0+p_1, 
$$
$$
 d(\underline{\mathbb S_{(P_3,0)}})=1-\beta_1=1-p_0,
$$
as well as
$$
\omega_1=p_1\frac{p_0+p_1}{p_1+p_2}, ~~~~~~~~~~~~~~~~~\omega_2=p_1\frac{p_1+p_2}{p_0+p_1}, ~~~~~~~~~~~~~~~~~\omega_3=p^2 _2, ~~~~~~~~~~~~~~~~~\omega_4=p^2 _0.
$$
In addition, for any step $k\in\mathbb N$, the following relationships hold: 
$$
N_{1,c_i}+N_{2,c_i}+N_{3,c_i}+N_{4,c_i}=2^{k}
$$
and
$$
N_{1,c_i}=N_{2,c_i}=N_{3,c_i}=N_{4,c_i}=2^{k-2}.
$$

The local structure our set can be characterized by the following scheme: 
$$
\begin{array}{ccccccc}
&{I_0}&  \\
\swarrow & &\searrow \\
\omega_2&  &  \omega_4 &
\end{array}
\begin{array}{ccccccc}
&{\omega_1}&  \\
\swarrow & &\searrow \\
\omega_2&  &  \omega_4 &
\end{array}
\begin{array}{ccccccc}
&{\omega_2}&  \\
\swarrow & &\searrow \\
\omega_1&  &  \omega_3 &
\end{array}
\begin{array}{ccccccc}
&{\omega_3}&  \\
\swarrow & &\searrow \\
\omega_1&  &  \omega_3 &
\end{array}
\begin{array}{ccccccc}
&{\omega_4}&  \\
\swarrow & &\searrow \\
\omega_2&  &  \omega_4 &
\end{array}
$$
Here $I_0=[\inf \mathbb S_{(-P_3,0)}, \mathbb S_{(-P_3,0)}]$.

So, the Hausdorff dimension $\dim_H(\mathbb S_{(-P_3,0)})$ of $\mathbb S_{(-P_3,0)}$ is equal to 
$$
\alpha_{*}=\liminf_{k\to\infty}{\gamma_k},
$$
where $(\gamma_k)$ is a sequence of numbers satisfying the equation
\begin{equation*}
\label{eq: the main for examples}
\left((\omega_{2})^{\gamma_1} +\left(\omega_{4}\right)^{\gamma_1} \right)\prod^k _{i=2}{\left({2^{i-2}}\left(\left(\omega_{1}\right)^{\gamma_i} + \left(\omega_{2}\right)^{\gamma_i}+\left(\omega_{3}\right)^{\gamma_i}+\left(\omega_{4}\right)^{\gamma_i}\right)  \right)}=1.
\end{equation*}

\begin{example} Suppose $p_0=\frac{1}{6}$, $p_1=\frac{2}{6}=\frac{1}{3}$, and $p_2=\frac{3}{6}=\frac{1}{2}$. Then the set 
$\mathbb S_{(P_3,0)}$ is a self-similar fractal whose Hausdorff dimension  is approximately equal to $0.408985$; but the set  $S_{(-P_3,0)}$ is not a self-similar fractal, its Hausdorff dimension  is approximately equal to $0.422592$.
\end{example}

\begin{example}
Suppose $p_0=p_2=0.25$ and $p_1=0.5$. Then the sets $\mathbb S_{(P_3,0)}$ and $\mathbb S_{(-P_3,0)}$ are self-similar fractals and their Hausdorff dimensions  approximately equal  $0.46496$.
\end{example}

\end{document}